\documentclass[10pt,twocolumn,letterpaper]{article}

\usepackage[pagenumbers]{cvpr} % To force page numbers, e.g. for an arXiv version

\usepackage{graphicx}
\usepackage{amsmath,amsthm}
\usepackage{amssymb}
\usepackage{booktabs}
\usepackage{xcolor,soul,verbatim}
\newcommand{\cbox}[2]{\colorbox{#1}{$\displaystyle #2$}}

\usepackage[pagebackref,breaklinks,colorlinks]{hyperref}

\newtheorem{thm}{Theorem}[section]
\newtheorem{dfn}[thm]{Definition}
\newtheorem{pro}[thm]{Problem}

\newtheorem{lem}[thm]{Lemma}
\newtheorem{exa}[thm]{Example}

\newcommand{\R}{\mathbb R}

\newcommand{\la}{\lambda}

\newcommand{\si}{\sigma}

\newcommand{\ep}{\varepsilon}

\newcommand{\HKS}{\mathrm{HKS}}

\newcommand{\ADD}{\mathrm{ADD}}
\newcommand{\ASD}{\mathrm{ASD}}
\newcommand{\SDM}{\mathrm{SDM}}

\newcommand{\SDV}{\mathrm{SDV}}

\newcommand{\VID}{\mathrm{VID}}
\newcommand{\VSM}{\mathrm{VSM}}

\newcommand{\LAC}{\mathrm{LAC}}
\newcommand{\PDD}{\mathrm{PDD}}

\newcommand{\SDD}{\mathrm{SDD}}
\newcommand{\MSD}{\mathrm{MSD}}
\newcommand{\WSD}{\mathrm{WSD}}
\newcommand{\WDD}{\mathrm{WDD}}
\newcommand{\wSDD}{\widetilde{\SDD}}
\newcommand{\RDD}{\mathrm{RDD}}
\newcommand{\wRDD}{\widetilde{\RDD}}

\newcommand{\AMD}{\mathrm{AMD}}

\newcommand{\EMD}{\mathrm{EMD}}

\newcommand{\lra}{\leftrightarrow}
\newcommand{\Lra}{\Leftrightarrow}
\newcommand{\es}{\emptyset}
\newcommand{\vl}{\,:\,}

\newcommand{\vect}[2]{ \left( \begin{array}{c} 
 #1 \\ #2 \end{array} \right)}

\newcommand{\mat}[4]{ \left( \begin{array}{cc} 
 #1 & #2 \\ #3 & #4 \end{array} \right)}
\newcommand{\matv}[4]{ \left( \begin{array}{cc} 
 #1 & #3 \\ #2 & #4 \end{array} \right)}

\usepackage[capitalize]{cleveref}
\crefname{section}{Sec.}{Secs.}
\Crefname{section}{Section}{Sections}
\Crefname{table}{Table}{Tables}
\crefname{table}{Tab.}{Tabs.}

\def\cvprPaperID{11825} % *** Enter the CVPR Paper ID here
\def\confName{CVPR}
\def\confYear{2023}

\begin{document}

%%%%%%%%% TITLE - PLEASE UPDATE
\title{Simplexwise Distance Distributions for finite spaces with metrics and measures}

\author{
Vitaliy Kurlin\\
Computer Science department\\
University of Liverpool, UK\\
{\tt\small vitaliy.kurlin@gmail.com}
}
\maketitle

%%%%%%%%% ABSTRACT
\begin{abstract}
A finite set of unlabelled points in Euclidean space is the simplest representation of many real objects from mineral rocks to sculptures.
Since most solid objects are rigid, their natural equivalence is rigid motion or isometry maintaining all inter-point distances.
More generally, any finite metric space is an example of a metric-measure space that has a probability measure and a metric satisfying all axioms.
\smallskip

This paper develops Simplexwise Distance Distributions (SDDs) for any finite metric spaces and metric-measures spaces.
These SDDs classify all known non-equivalent spaces that were impossible to distinguish by simpler invariants.
We define metrics on SDDs that are Lipschitz continuous and allow exact computations whose parametrised complexities are polynomial in the number of given points.  
\end{abstract}

%1==========================
\section{Motivations for classifying metric spaces}
\label{sec:intro}

The simplest representation of any rigid object such a car or a sculpture is a finite set (\emph{cloud}) $S\subset\R^n$ of $m$ \emph{unlabelled} points, where $n=2,3$ are the most practical dimensions. 
\smallskip

The rigidity of many solid objects motivates us to study them up to \emph{rigid motion}, which is a composition of translations and rotations in Euclidean space $\R^n$.
We can consider any finite set $X$ with a \emph{metric} that is a distance function $d_X:X\times X\to[0,+\infty)$ satisfying all metric axioms.
The natural equivalence relation on metric spaces is an \emph{isometry} that is any map $f:X\to Y$ maintaining all inter-point distances so that $d_X(p,q)=d_Y(f(p),f(q))$ for $p,q\in X$.
\smallskip

In $\R^n$, any isometry is a composition of a mirror reflection with some rigid motion.
Any orientation-preserving isometry can be realised as a continuous rigid motion.
\smallskip

The \emph{shape} of a rigid object is mathematically defined as its isometry class.
Any non-rigid deformation defines a weaker equivalence (than isometry) with a smaller space of flexible shapes.
Comparing shapes up to isometry requires finer invariants to distinguish many more isometry classes.
\smallskip

The mathematical approach to distinguish spaces up to isometry uses \emph{invariants} that are properties preserved by any isometry.
Any invariant $I$ maps all isometric spaces to the same value, hence has \emph{no false negatives} that are pairs of isometric spaces $S\cong Q$ with $I(S)\neq I(Q)$.
\smallskip

A \emph{complete} invariant $I$ should distinguish all non-isometric clouds, so if $S\not\cong Q$ then $I(S)\neq I(Q)$.
Equivalently, $I$ has \emph{no false positives} that are pairs of non-isometric spaces $S\not\cong Q$ with $I(S)=I(Q)$.
Then $I$ is a DNA-style code that uniquely identifies any space $S$ up to isometry. 
\smallskip

Since real data are always noisy, a useful complete invariant must be also continuous under the movement of points. 
Satisfying both completeness and continuity is extremely challenging for sets of $m$ \emph{unlabelled} points because of $m!$ potential permutations that match all $m$ points. 
\smallskip

A complete and continuous invariant for $m=3$ points consists of three pairwise distances $a,b,c$ (sides of a triangle) and is known in school as the SSS theorem \cite{weisstein2003triangle} about the congruence (isometry) of triangles.
As a result, the isometry space of 3-point sets is continuously mapped as a quadrangular cone $\{0<a\leq b\leq c\leq a+b\}$ parametrised by $a\leq b\leq c$ satisfying one triangle inequality $a+b\geq c$.
\smallskip

\begin{figure}[h!]
\centering
\includegraphics[width=\linewidth]{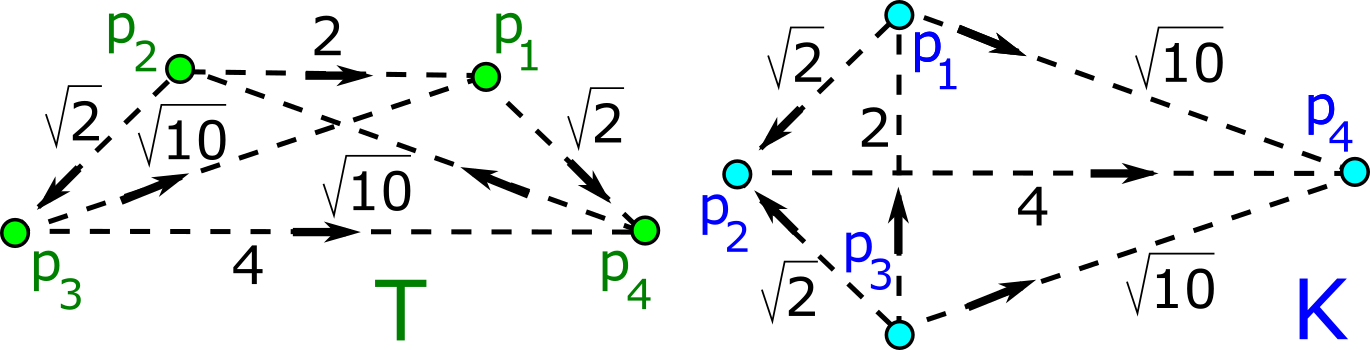}
\caption{
\textbf{Left}: trapezium $T=\{(1,1),(-1,1),(-2,0),(2,0)\}$. 
\textbf{Right}: kite $K=\{(0,1),(-1,0),(0,-1),(3,0)\}$. 
Both $T$ and $K$ have the same six pairwise distances $\sqrt{2},\sqrt{2},2,\sqrt{10},\sqrt{10},4$.
}
\label{fig:4-point_clouds}
\end{figure}

The full description above had no easy analogue for $m\geq 4$ points in $\R^n$.
One obstacle was a family of 4-point sets in $\R^2$ that have the same six pairwise distances, see Fig.~\ref{fig:4-point_clouds}.

\begin{pro}[complete isometry invariants with computable continuous metrics]
\label{pro:isometry}
Design an invariant $I$ of finite metric spaces satisfying the following properties:
\medskip

\noindent
\textbf{(a)}
\emph{completeness} : 
$S,Q$ are isometric $\Lra$ $I(S)=I(Q)$; 
\medskip

\noindent
\textbf{(b)}
\emph{Lipschitz continuity} :
if any point of $C$ is perturbed within its $\ep$-neighbourhood then
 $I(S)$ changes by at most $\la\ep$ for a constant $\la$ and a metric $d$ 
 satisfying all axioms:
\smallskip

\noindent 
(1) $d(I(S),I(Q))=0$ if and only if $S\cong Q$ are isometric,
\smallskip

\noindent 
(2) \emph{symmetry} : $d(I(S),I(Q))=d(I(Q),I(S))$,
\smallskip

\noindent 
%\emph{triangle inequality} : 
(3) $d(I(S),I(Q))+d(I(Q),I(T))\geq d(I(S),I(T))$;
\medskip

\noindent
\textbf{(c)}
\emph{computability} : $I(S)$ and $d$ are computable in a polynomial time in the number $m$ of points in given spaces.
%\bs
\end{pro}

Condition~(\ref{pro:isometry}b) asking for a continuous metric is stronger than the completeness in (\ref{pro:isometry}a).
Detecting an isometry $S\cong Q$ gives a discontinuous metric, say $d=1$ for all non-isometric $S\not\cong Q$ even if $S,Q$ are nearly identical.
Any metric $d$ satisfying the first axiom in (\ref{pro:isometry}b) detects an isometry $S\cong Q$ by checking if $d(I(S),I(Q))=0$.
\smallskip

Problem~\ref{pro:isometry} was open at least since 1974 when Gilbert and Shepp  \cite{gilbert1974textures} described a 4-parameter family of 4-point sets in $\R^2$ that have the same six pairwise distances. 
Problem~\ref{pro:isometry} is also motivated by the weaknesses \cite{smith2022families,elkin2020mergegram,elkin2021isometry} of persistent homology in Topological Data Analysis.
\smallskip

Section~\ref{sec:past} reviews the closely related work on invariants of point sets and more general metric spaces.
Section~\ref{sec:SDD} introduces the Simplexwise Distance Distributions ($\SDD$s), which substantially generalise all past distance-based invariants of finite metric spaces. 
Section~\ref{sec:examples} shows that $\SDD$s are simple enough for manual computations and classifying infinite families of clouds that cannot be distinguished by simpler distance distributions.   
Section~\ref{sec:metrics} develops Lipschitz continuous metrics on $\SDD$s that are computed in a parametrized polynomial time in the number $m$ of points.
\smallskip

We consider Problem~\ref{pro:isometry} a first important step toward understanding moduli spaces of any data objects.
Metric spaces and isometry can be replaced by other data and equivalence, respectively, to get analogues of Problem~\ref{pro:isometry}.
This paper extends section 3 of \cite{widdowson2023recognizing}, whose 8-page version without proofs and big examples will appear soon.
In the papers \cite{widdowson2021pointwise,widdowson2022resolving,widdowson2023recognizing} the first author implemented all algorithms, the second author designed all theory, proofs, and examples. 

%2==========================
\section{Past work on isometries and metric spaces}
\label{sec:past}

This section reviews the related work starting from a simpler version of Problem~\ref{pro:isometry} asking only to detect a potential isometry between clouds of $m$ unlabelled points
\smallskip

\noindent  
\textbf{Isometry detection} refers to a simpler version of Problem~\ref{pro:isometry} to algorithmically detect a potential isometry between given clouds of $m$ points in $\R^n$. 
%The first solution by Alt et al \cite{alt1988congruence} required $O(m^{n-2}\log m)$ time.
The best algorithm by Brass and Knauer \cite{brass2000testing} takes $O(m^{\lceil n/3\rceil}\log m)$ time, so $O(m\log m)$ in $\R^3$ \cite{brass2004testing}.
The latest advance is the $O(m\log m)$ algorithm in $\R^4$ \cite{kim2016congruence}.
These algorithms output a binary answer (yes/no) without quantifying similarity between non-isometric clouds by a continuous metric.
\smallskip

\noindent
\textbf{Multidimensional scaling} (MDS).
For a given $m\times m$ distance matrix of any $m$-point cloud $A$, MDS \cite{schoenberg1935remarks} finds an embedding $A\subset\R^k$ (if it exists) preserving all distances of $M$ for a dimension $k\leq m$.
%The $m$ eigenvalues of the Gram matrix expressed via $D$ need $O(m^3)$ time.
A final embedding $A\subset\R^k$ uses %orthonormal 
eigenvectors whose ambiguity up to signs gives an exponential comparison time that can be close to $O(2^m)$.
\smallskip

\noindent
\textbf{The Heat Kernel Signature} ($\HKS$) 
%$\HKS_M:\R^+\times M\to\R$ 
is a complete isometry invariant of a manifold $M$ whose the Laplace-Beltrami operator has distinct eigenvalues by \cite[Theorem~1]{sun2009concise}.
%This completeness reduces the detection of an isometry  $M\to N$ to the test if $\HKS_M(x,t)=\HKS_N(f(x),t)$ for a suitable bi-continuous bijection (homeomorphism) $f:M\to N$ and all (infinitely many values of) $x\in M$ and $t>0$. 
If $M$ is sampled by points, $\HKS$ can be discretized and remains continuous \cite[section~4]{sun2009concise} but the completeness is unclear.
\smallskip

\noindent
\textbf{The Hausdorff distance} \cite{hausdorff1919dimension} can be defined for any subsets $A,B$ in an ambient metric space as  $d_H(A,B)=\max\{\vec d_H(A,B), \vec d_H(B,A) \}$, where the directed Hausdorff distance is 
$\vec d_H(A,B)=\sup\limits_{p\in A}\inf\limits_{q\in B}|p-q|$.
To take into account isometries, one can minimize the Hausdorff distance over all isometries \cite{huttenlocher1993comparing,chew1992improvements,chew1999geometric}.
 % $f$ from the full Euclidean group $\iso(\R^n)$.
%For $n=1$, the Hausdorff distance minimized over translations in $\R$ for sets of at most $m$ points can be found in time $O(m\log m)$ \cite{rote1991computing}.
For $n=2$, the Hausdorff distance minimized over isometries in $\R^2$ for sets of at most $m$ point needs $O(m^5\log m)$ time \cite{chew1997geometric}. 
For a given $\ep>0$ and $n>2$, the related problem to decide if $d_H\leq\ep$ up to translations has the time complexity $O(m^{\lceil(n+1)/2\rceil})$ \cite[Chapter~4, Corollary~6]{wenk2003shape}. 
For general isometry, only approximate algorithms tackled minimizations for infinitely many rotations initially in $\R^3$ \cite{goodrich1999approximate} and in $\R^n$ \cite[Lemma~5.5]{anosova2022algorithms}. 
\smallskip

\noindent
\textbf{The Gromov-Wasserstein distances} can be defined for metric-measure spaces, not necessarily sitting in a common ambient space.
The simplest Gromov-Hausdorff (GH) distance cannot be approximated with any factor less than 3 in polynomial time unless P = NP \cite[Corollary~3.8]{schmiedl2017computational}.
Polynomial-time algorithms for GH were designed for ultrametric spaces  \cite{memoli2021gromov}.
However, GH spaces are challenging even for finite point sets in the line $\R$, see \cite{majhi2019approximating} and \cite{zava2023gromov}.
\smallskip

\noindent  
\textbf{Experimental approaches} cover a wide variety of descriptors designed manually or optimised through machine learning, for example, Scale Invariant Feature Transform \cite{toews2013efficient,rister2017volumetric,spezialetti2019learning,zhu2022point}.
Some of these descriptors are designed for invariance under permutations of points \cite{qi2017pointnet++,zaheer2017deep}, and also for invariance under isometry \cite{chen2021equivariant,nigam2022equivariant,simeonov2022neural}, for example, in Geometric Deep Learning \cite{bronstein2017geometric,bronstein2021geometric}.
Among many obstacles \cite{dong2018boosting,akhtar2018threat,laidlaw2019functional,guo2019simple,colbrook2022difficulty}, the hard one is to theoretically guarantee the completeness and Lipschitz continuity of such descriptors under perturbations as in Problem~\ref{pro:isometry}.
\smallskip

\noindent
\textbf{Local distributions of distances} in M\'emoli's seminal work \cite{memoli2011gromov,memoli2022distance} for metric-measure spaces, or shape distributions \cite{osada2002shape, belongie2002shape, grigorescu2003distance, manay2006integral, pottmann2009integral}, are first-order versions of the new $\SDD$s.
Another approach to Problem~\ref{pro:isometry} uses direction-based invariants \cite{kurlin2022computable}, which inspired Complete Neural Networks \cite{hordan2023complete}.
The Lipschitz continuity was proved \cite[Theorem~4.9]{kurlin2022computable} in general position but not for near-singular configurations, for example, when a triangle degenerates to a line.
These degeneracies will be addressed in the forthcoming work \cite{kurlin2023strength} extending $\SDD$ to a complete invariant of clouds in $\R^n$.  

%3==========================
\section{Simplexwise Distance Distribution (SDD)}
\label{sec:SDD}

This section introduces the isometry invariants $\SDD$ for a finite cloud $C$ of unlabelled points in any metric space $X$.  
The \emph{lexicographic} order $u<v$ on vectors $u=(u_1,\dots,u_h)$ and $v=(v_1,\dots,v_h)$ means that if the first $i$ (possibly, $i=0$) coordinates of $u,v$ coincide then $u_{i+1}<v_{i+1}$.
%For example, $(1,2)<(2,1)<(2,2)$.
Let $S_h$ denote the permutation group on indices $1,\dots,h$.

\begin{dfn}[$\RDD(C;A)$] % for a subset $A\subset C$]
\label{dfn:RDD}
Let $C$ be a cloud of $m$ unlabelled points in a space with a metric $d$.
A \emph{basis} sequence $A=(p_1,\dots,p_h)\in C^h$ consists of $1\leq h<m$ distinct points.
Let $D(A)$ be the triangular \emph{distance} matrix whose entry $D(A)_{i,j-1}$ is $d(p_i,p_j)$ for $1\leq i<j\leq h$, all other entries are filled by zeros.
Any permutation $\xi\in S_h$ acts on $D(A)$ by mapping $D(A)_{ij}$  to $D(A)_{kl}$, where $k\leq l$ is the pair of indices $\xi(i),\xi(j)-1$ written in increasing order.
\smallskip

For any other point $q\in C-A$, write distances from $q$ to $p_1,\dots,p_h$ as a column.
The $h\times (m-h)$-matrix $R(C;A)$ is formed by these $m-h$ lexicographically ordered columns.
The action of $\xi$ on $R(C;A)$ maps any $i$-th row to the $\xi(i)$-th row, after which all columns can be written again in the lexicographic order.
The \emph{Relative Distance Distribution} $\RDD(C;A)$ is the equivalence class of the pair $[D(A),R(C;A)]$ of matrices up to  permutations $\xi\in S_h$.
%\bs
\end{dfn}

For $h=1$ and a basis sequence $A=(p_1)$, the matrix $D(A)$ is empty and $R(C;A)$ is a single row of distances (in the increasing order) from $p_1$ to all other points $q\in C$.
For $h=2$ and a basis sequence $A=(p_1,p_2)$, the matrix $D(A)$ is the single number $d(p_1,p_2)$ and $R(C;A)$ consists of two rows of distances from $p_1,p_2$ to all other points $q\in C$.

\begin{figure}[h!]
\centering
\includegraphics[width=\linewidth]{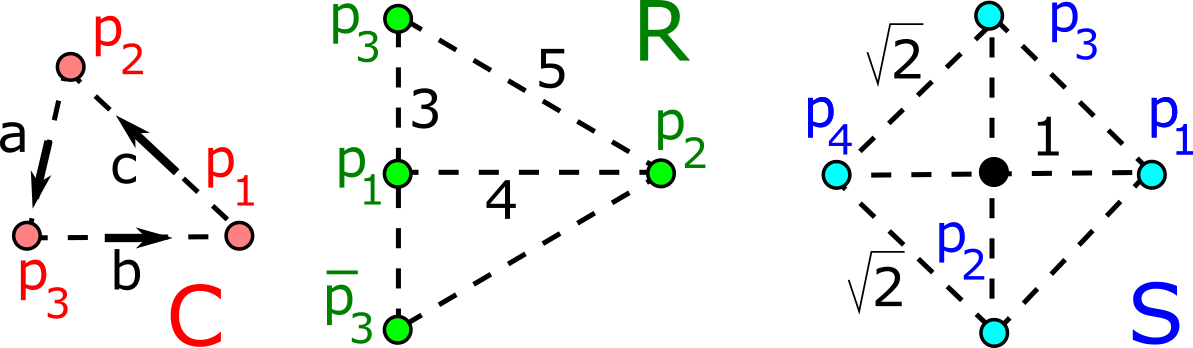}
\caption{\textbf{Left}: a cloud $C=\{p_1,p_2,p_3\}$ with distances $a\leq b\leq c$.
\textbf{Middle}: the triangular cloud $R=\{(0,0),(4,0),(0,3)\}$. 
\textbf{Right}: the square cloud $S=\{(1,0),(-1,0),(0,1),(-1,0)\}$. 
}
\label{fig:triangular_clouds}
\end{figure}

\begin{exa}[$\RDD$ for a 3-point cloud $C$]
\label{exa:RDD}
Let $C\subset\R^2$ consist of $p_1,p_2,p_3$ with inter-point distances $a\leq b\leq c$ ordered counter-clockwise as in Fig.~\ref{fig:triangular_clouds}~(left).  
Then
$$\RDD(C;p_1)=[\es;(b,c)], 
\RDD(C;\vect{p_2}{p_3})=[a;\vect{c}{b}],$$
$$\RDD(C;p_2)=[\es;(a,c)], 
\RDD(C;\vect{p_3}{p_1})=[b;\vect{a}{c}],$$
$$\RDD(C;p_3)=[\es;(a,b)], 
\RDD(C;\vect{p_1}{p_2})=[c;\vect{b}{a}].$$
We have written $\RDD(C;A)$ for a basis sequence $A=(p_i,p_j)$ of ordered points represented by a column.
Swapping the points $p_1\lra p_2$ makes the last $\RDD$ above equivalent to 
$\RDD\big(C;\vect{p_2}{p_1}\big)=\big[c;\vect{a}{b}\big]$.
%\bs
\end{exa}

Though $\RDD(C;A)$ is defined up to a permutation $\xi\in S_h$ of $h$ points in $A\in C^h$, comparisons of $\RDD$s will be practical for $h=2,3$ with metrics independent of $\xi$. 

\begin{dfn}[Simplexwise Distance Distribution $\SDD(C;h)$]
\label{dfn:SDD}
Let $C$ be a cloud of $m$ unlabelled points in a metric space.
For an integer $1\leq h<m$, the \emph{Simplexwise Distance Distribution} $\SDD(C;h)$ is the unordered set of $\RDD(C;A)$ for all unordered $h$-point subsets $A\subset C$.
%\bs
\end{dfn}

For $h=1$ and any $m$-point cloud $C$, the distribution $\SDD(C;1)$ can be considered as a matrix of $m$ rows of ordered distances from every point $p\in C$ to all other $m-1$ points.
If we lexicographically order these $m$ rows and collapse any $l>1$ identical rows into a single one with the weight $l/m$, then we get the
Pointwise Distance Distribution $\PDD(C;m-1)$ introduced in \cite[Definition~3.1]{widdowson2022resolving}. 
\smallskip

The PDD was simplified to the easier-to-compare vector of Average Minimum Distances \cite{widdowson2022average}:
$\AMD_k(C)=\dfrac{1}{m}\sum\limits_{i=1}^m d_{ik}$, where $d_{ik}$ is the distance from a point $p_i\in C$ to its $k$-th nearest neighbor in $C$.
These neighbor-based invariants can be computed in a near-linear time in $m$ \cite{elkin2022new} and were pairwise compared for all all 660K+ periodic crystals in the world's largest database of real materials \cite{widdowson2022resolving}. 
Definition~\ref{dfn:SDM} similarly maps $\SDD$ to a smaller invariant.
\smallskip

Recall that the 1st moment of a set of numbers $a_1,\dots,a_k$ is the \emph{average} $\mu=\dfrac{1}{k}\sum\limits_{i=1}^k a_i$.
The 2nd moment is the \emph{standard deviation} $\si=\sqrt{\dfrac{1}{k}\sum\limits_{i=1}^k (a_i-\mu)^2}$.
For $l\geq 3$, the $l$-th \emph{standardized moment} \cite[section~2.7]{keeping1995introduction} is $\dfrac{1}{k}\sum\limits_{i=1}^k \left(\dfrac{a_i-\mu}{\si}\right)^l$.
 
\begin{dfn}[Simplexwise Distance Moments $\SDM$]
\label{dfn:SDM}
For any $m$-point cloud $C$ in a metric space,
 let $A\subset C$ be a subset of $h$ unordered points.
The \emph{Sorted Distance Vector} $\SDV(A)$ is the list of all $\frac{h(h-1)}{2}$ pairwise distances between points of $A$ written in increasing order. 
The vector $\vec R(C;A)\in\R^{m-h}$ is obtained from the $h\times(m-h)$ matrix $R(C;A)$ in Definition~\ref{dfn:RDD} by writing the vector of $m-h$ column averages in increasing order.
\smallskip

The pair $[\SDV(A);\vec R(C;A)]$ is the \emph{Average Distance Distribution} $\ADD(C;A)$ considered as a vector of length $\frac{h(h-3)}{2}+m$.
The unordered collection of $\ADD(C;A)$ for all $\binom{m}{h}$ unordered subsets $A\subset C$ is the Average Simplexwise Distribution $\ASD(C;h)$.
The \emph{Simplexwise Distance Moment} $\SDM(C;h,l)$ is the $l$-th (standardized for $l\geq 3$) moment of $\ASD(C;h)$ considered as a probability distribution of $\binom{m}{h}$ vectors, separately for each coordinate.
%\bs
\end{dfn}

\begin{exa}[$\SDD$ and $\SDM$ for $T,K$]
\label{exa:SDD}
Fig.~\ref{fig:4-point_clouds} shows the non-isometric 4-point clouds $T,K$ with the same Ordered Pairwise Distances: $\SDV=\{\sqrt{2},\sqrt{2},2,\sqrt{10},\sqrt{10},4\}$, see infinitely many examples in \cite{boutin2004reconstructing}.
The arrows on the edges of $T,K$ show orders of points in each pair of vertices for $\RDD$s.
Then $T,K$ are distinguished up to isometry by $\SDD(T;2)\neq \SDD(K;2)$ in Table~\ref{tab:SDD+TK}.
The 1st coordinate of $\SDM(C;2,1)\in\R^3$ is the average of the six distances from $\SDV$ (the same for $T,K$) but the other two coordinates (column averages from $R(C;A)$ matrices) differ.
%\bs
\end{exa}

\begin{table}
  \centering
  \begin{tabular}{@{}l|l@{}}
    \toprule
   	$\RDD(T;A)$ in $\SDD(T;2)$ & $\RDD(K;A)$ in $\SDD(K;2)$ \\
    
    	$[\sqrt{2},\matv{2}{\cbox{yellow}{\sqrt{10}}}{\sqrt{10}}{4}]\times 2$ &
	$[\sqrt{2},\matv{2}{\cbox{yellow}{\sqrt{2}}}{\sqrt{10}}{4}]\times 2$ 
	\smallskip \\
	
$[2,\matv{\sqrt{2}}{\cbox{yellow}{\sqrt{10}}}{\sqrt{10}}{\cbox{yellow}{\sqrt{2}}}]$ & 
$[2,\matv{\sqrt{2}}{\cbox{yellow}{\sqrt{2}}}{\sqrt{10}}{\cbox{yellow}{\sqrt{10}}}]$ 
\smallskip \\
    
	$[\sqrt{10},\matv{\sqrt{2}}{\cbox{yellow}{2}}{\cbox{yellow}{4}}{\cbox{yellow}{\sqrt{2}}}]\times 2$ &
	$[\sqrt{10},\matv{\sqrt{2}}{\cbox{yellow}{4}}{\cbox{yellow}{2}}{\cbox{yellow}{\sqrt{10}}}]\times 2$ \smallskip \\
	
	    $[4,\matv{\sqrt{2}}{\sqrt{10}}{\cbox{yellow}{\sqrt{10}}}{\cbox{yellow}{\sqrt{2}}}]$ &
	$[4,\matv{\sqrt{2}}{\sqrt{10}}{\cbox{yellow}{\sqrt{2}}}{\cbox{yellow}{\sqrt{10}}}]$  \\

    \midrule
  	$\ADD(T;A)$ in $\ASD(T;2)$ & $\ADD(K;A)$ in $\ASD(K;2)$ \\
    $[\sqrt{2},(\frac{2+\cbox{yellow}{\sqrt{10}}}{2},\frac{4+\sqrt{10}}{2})]\times 2$ &
    $[\sqrt{2},(\frac{2+\cbox{yellow}{\sqrt{2}}}{2},\frac{4+\sqrt{10}}{2})]\times 2$ \\

    $[2,(\cbox{yellow}{\frac{\sqrt{2}+\sqrt{10}}{2},\frac{\sqrt{2}+\sqrt{10}}{2}})]$ &
    $[2,(\cbox{yellow}{\sqrt{2},\sqrt{10}})]$ \\

    $[\sqrt{10},(\frac{2+\cbox{yellow}{\sqrt{2}}}{2},\frac{4+\sqrt{2}}{2})]\times 2$ &
    $[\sqrt{10},(\frac{2+\cbox{yellow}{\sqrt{10}}}{2},\frac{4+\sqrt{2}}{2})]\times 2$ \\
    
       $[4,(\frac{\sqrt{2}+\sqrt{10}}{2},\frac{\sqrt{2}+\sqrt{10}}{2})]$ & $[4,(\frac{\sqrt{2}+\sqrt{10}}{2},\frac{\sqrt{2}+\sqrt{10}}{2})]$     \\

    \midrule
   	$\SDM_1=\dfrac{3+\sqrt{2}+\sqrt{10}}{3}$ & 
   	$\SDM_1=\dfrac{3+\sqrt{2}+\sqrt{10}}{3}$ \\

   	$\SDM_2=\dfrac{\cbox{yellow}{6+2\sqrt{2}+4\sqrt{10}}}{12}$ & 
   	$\SDM_2=\dfrac{\cbox{yellow}{8+5\sqrt{2}+3\sqrt{10}}}{12}$ \\

   	$\SDM_3=\frac{16+\cbox{yellow}{4\sqrt{2}+4\sqrt{10}}}{12}$ & 
   	$\SDM_3=\frac{16+\cbox{yellow}{3\sqrt{2}+5\sqrt{10}}}{12}$ \\

    \bottomrule
  \end{tabular}
  \caption{The Simplexwise Distance Distributions from Definition~\ref{dfn:SDD} for the 4-point clouds $T,K\subset\R^2$ in Fig.~\ref{fig:4-point_clouds}. 
The symbol $\times 2$ indicates a doubled $\RDD$.
The three bottom rows show coordinates of $\SDM(C;2,1)\in\R^3$ from Definition~\ref{dfn:SDM} for $h=2$, $l=1$ and both $C=T,K$.
Different elements are \hl{highlighted}.
}
\label{tab:SDD+TK}
\end{table}

Some of the $\binom{m}{h}$ $\RDD$s in $\SDD(C;h)$ can concide as in Example~\ref{exa:SDD}.
If we collapse any $l>1$ identical $\RDD$s into a single $\RDD$ with the \emph{weight} $l/\binom{m}{h}$, $\SDD$ can be considered as a weighted probability distribution 
of $\RDD$s.
\smallskip

All time complexities are proved for a random-access machine (RAM) model.
In a general metric space, a point cloud $C$ is usually given by a distance matrix on (arbitrarily ordered) points of $C$. 
Hence we assume that the distance between any points of $C$ is accessible in a constant time.

\begin{thm}[invariance and time of $\SDD$s]
\label{thm:SDD_time}
For $h\geq 1$ and any cloud $C$ of $m$ unlabelled points in a metric space, $\SDD(C;h)$ is an isometry invariant, which can be computed in time $O(m^{h+1}/(h-1)!)$.
For any $l\geq 1$, the invariant $\SDM(C;h,l)\in\R^{m+\frac{h(h-3)}{2}}$ has the same time.
\end{thm}
\begin{proof}
Any isometry $S\to Q$ preserves distances, hence induces a bijection 
$\SDD(S;h)\to\SDD(Q;h)$ for $h\geq 1$.
\smallskip

By Definition~\ref{dfn:SDD}, for any $h\geq 1$ and a cloud $C$ of $m$ unlabelled points in a metric space, the Simplexwise Distance Distribution $\SDD(C;h)$ of consists of $\binom{m}{h}=\frac{m!}{h!(m-h)!}=O(m^h/h!)$ Relative Distance Distributions $\RDD(C;A)$ for any unordered subset $A\subset C$ of $h$ points.
\smallskip

For any order of points of $A$, every $\RDD(C;A)$ consists of the distance matrix $D(A)$, which needs $O(h^2)$ time and $h\times(m-h)$ matrix $R(C;A)$, which needs $h(m-h)$ time.
Since $h\leq m$, the extra factor $O(hm)$ gives the final time $O(m^{h+1}/(h-1)!)$ for $\SDD(C;h)$.
\smallskip

For a fixed $h$-point subset $A\subset C$, the vector $\vec R(C;A)$ from Definition~\ref{dfn:SDM} needs $O(hm)$ time to average $h$ distances in $m-h$ columns and $O(m\log m)$ time to order these averages.
The list $\SDV(A)$ of Ordered Pairwise Distances is obtained by sorting all pairwise distances from $D(A)$ in time $O(h^2\log h)$.
So the Average Distance Distribution $\ADD(C;A)$ obtained by concatenating the ordered vectors $\SDV(A)\in\R^{\frac{h(h-1)}{2}}$ and $\vec R(C;A)\in\R^{m-h}$ requires only $O((h^2+m)\log m)$ extra time.
%For $1\leq h<m$, we roughly estimate $(h^2+m)\log m$ as $O() 
Hence the Average Simplexwise Distribution $\ASD(C;h)$ for all $h$-point subsets $A\subset C$ needs $O(m^{h+1}/(h-1)!)$ time including $O((h^2+m)\log m)$, the same as the initial $\SDD(C;h)$.
\smallskip

For $l=1$, the first raw moment $\SDM(C;h,1)$ is the simple average of all $k=\binom{m}{h}$ vectors $\ADD(C;A)$ of length $O(hm)$, hence needs $O(m^{h+1}/(h-1)!)$ time.
For $l=2$, the standard deviation $\si$ of each coordinate in all vectors $\ADD(C;A)$ requires the same time.
Then, for any fixed $l\geq 3$, the $l$-th standardized moment  
 $\dfrac{1}{k}\sum\limits_{i=1}^k \left(\dfrac{a_i-\mu}{\si}\right)^l$
 needs again the same time $O(m^{h+1}/(h-1)!)$.
\end{proof}

We conjecture that $\SDD(C;h)$ is a complete isometry invariant of a cloud $C\subset\R^n$ for some $h\leq n$.
Section~\ref{sec:examples} shows that $\SDD(C;2)$ distinguishes all infinitely many known pairs \cite[Fig.~S4]{pozdnyakov2020incompleteness} of non-isometric $m$-point sets $S,Q\subset\R^3$ that have equal $\PDD(S)=\PDD(Q)$ 

\begin{figure*}[h!]
\centering
\includegraphics[width=\linewidth]{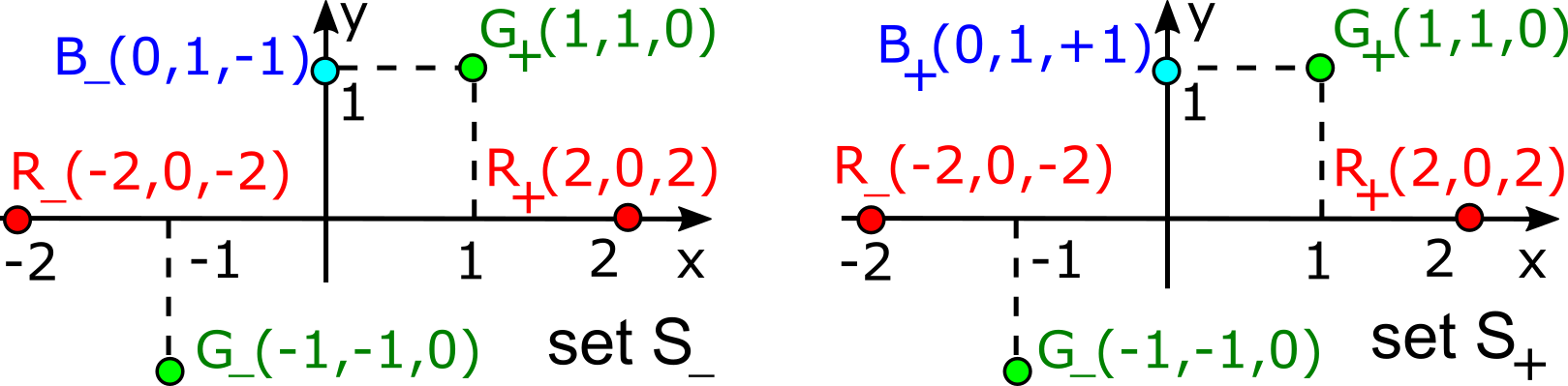}
\caption{\textbf{Left}: $(x,y)$-projection of the 5-point set $S_-\subset\R^3$ consisting of the green points $G_-=(-1,-1,0)$ and $G_{+}=(1,1,0)$, the red points $R_-=(-2,0,-2)$ and $R_+=(2,0,2)$, and the blue point $B_-=(0,1,-1)$.
\textbf{Right}: to get the set $S_+\subset\R^3$ from $S_-$, replace the point $B_-$ with the new point $B_+=(0,1,1)$.}
\label{fig:5-point_sets}
\end{figure*}

%4==========================
\section{The strength of the isometry invariant SDD}
\label{sec:examples}

Examples~\ref{exa:5-point_sets} and~\ref{exa:7-point_sets} distinguish clouds of 5 points and 7 points, respectively, in $\R^3$ by comparing their $\SDD$s of order 2.
In Example~\ref{exa:6-point_sets}, $\SDD(C;2)$ distinguishes 6-point clouds in a family of pairs depending on three parameters.  

%B1-----------------------------------------

\begin{exa}[5-point clouds]
\label{exa:5-point_sets}
Fig.~\ref{fig:5-point_sets} shows the 5-point clouds $S_{\pm}\subset\R^3$ taken from \cite[Figure S4(A)]{pozdnyakov2020incompleteness}.
\smallskip

\begin{table}[h!]
\begin{tabular}{@{}l|ccccc@{}}
 \toprule
distances of \hl{$S_-$} & $R_-$ & $R_+$ & $G_-$ & $G_+$ & \hl{$B_-$} \\
\midrule
$R_-(-2,0,-2)$ & $0$ & $\sqrt{32}$ & $\sqrt{6}$ & $\sqrt{14}$ & \hl{$\sqrt{6}$} \smallskip \\
$R_+(+2,0,+2)$ & $\sqrt{32}$ & 0 & $\sqrt{14}$ & $\sqrt{6}$ & \hl{$\sqrt{14}$} \smallskip \\
$G_-(-1,-1,0)$ & $\sqrt{6}$ & $\sqrt{14}$ & 0 & $\sqrt{8}$ & $\sqrt{6}$ \smallskip \\
$G_+(+1,+1,0)$ & $\sqrt{14}$ & $\sqrt{6}$ & $\sqrt{8}$ & 0 & $\sqrt{2}$ \smallskip \\
\hl{$B_-(0,+1,-1)$} & \hl{$\sqrt{6}$} & \hl{$\sqrt{14}$} & $\sqrt{6}$ & $\sqrt{2}$ & 0 \\
\bottomrule
\end{tabular}
\smallskip

\begin{tabular}{@{}l|ccccc@{}}
 \toprule
distances of $S_+$ & $R_-$ & $R_+$ & $G_-$ & $G_+$ & \hl{$B_+$} \\
\midrule
$R_-(-2,0,-2)$ & $0$ & $\sqrt{32}$ & $\sqrt{6}$ & $\sqrt{14}$ & \hl{$\sqrt{14}$} \smallskip \\
$R_+(+2,0,+2)$ & $\sqrt{32}$ & 0 & $\sqrt{14}$ & $\sqrt{6}$ & \hl{$\sqrt{6}$} \smallskip \\
$G_-(-1,-1,0)$ & $\sqrt{6}$ & $\sqrt{14}$ & 0 & $\sqrt{8}$ & $\sqrt{6}$ \smallskip \\
$G_+(+1,+1,0)$ & $\sqrt{14}$ & $\sqrt{6}$ & $\sqrt{8}$ & 0 & $\sqrt{2}$ \smallskip \\
\hl{$B_+(0,+1,+1)$} & \hl{$\sqrt{14}$} & \hl{$\sqrt{6}$} & $\sqrt{6}$ & $\sqrt{2}$ & 0 \\
\bottomrule
\end{tabular}
\caption{Distances between all points of the set $S_{\mp}$ in Fig.~\ref{fig:5-point_sets}.}
\label{tab:distances_S-+}
\end{table}

The sets $S_{\pm}$ are not isometric, because $S_+$ has the triple of points $ B_+, G_+, R_+$ with pairwise distances $\sqrt{2},\sqrt{6},\sqrt{6}$, but $S_-$ has no such a triple.
Table~\ref{tab:distances_S-+} \hl{highlights} differences between distance matrices.
If we order distances to neighbors, the matrices in Table~\ref{tab:ordered_distances_S-+} differ only in one pair.

\begin{table}[h]
\centering
\begin{tabular}{@{}l|cccc@{}}
\toprule
\hl{$S_-$} distances to & 1st neighbor & 2nd & 3rd & 4th \\
\midrule
$R_-=(-2,0,-2)$ & $\sqrt{6}$ & \hl{$\sqrt{6}$}  & $\sqrt{14}$ & $\sqrt{32}$ \smallskip \\
$R_+=(+2,0,+2)$ & $\sqrt{6}$ & \hl{$\sqrt{14}$}  & $\sqrt{14}$ & $\sqrt{32}$ \smallskip \\
$G_-=(-1,-1,0)$  & $\sqrt{6}$ & $\sqrt{6}$ & $\sqrt{8}$ & $\sqrt{14}$  \smallskip \\
$G_+=(+1,+1,0)$ & $\sqrt{2}$ & $\sqrt{6}$ & $\sqrt{8}$ & $\sqrt{14}$ \smallskip \\
$B_-=(0,+1,-1)$ & $\sqrt{2}$ & $\sqrt{6}$ & $\sqrt{6}$ & $\sqrt{14}$ \\
\bottomrule
\end{tabular}
\smallskip

\begin{tabular}{@{}l|cccc@{}}
\toprule
\hl{$S_+$} distances to & 1st neighbor & 2nd & 3rd & 4th \\
\midrule
$R_-=(-2,0,-2)$ & $\sqrt{6}$ & \hl{$\sqrt{14}$}  & $\sqrt{14}$ & $\sqrt{32}$ \smallskip \\
$R_+=(+2,0,+2)$ & $\sqrt{6}$ & \hl{$\sqrt{6}$}  & $\sqrt{14}$ & $\sqrt{32}$ \smallskip \\
$G_-=(-1,-1,0)$  & $\sqrt{6}$ & $\sqrt{6}$ & $\sqrt{8}$ & $\sqrt{14}$  \smallskip \\
$G_+=(+1,+1,0)$ & $\sqrt{2}$ & $\sqrt{6}$ & $\sqrt{8}$ & $\sqrt{14}$ \smallskip \\
$B_+=(0,+1,-1)$ & $\sqrt{2}$ & $\sqrt{6}$ & $\sqrt{6}$ & $\sqrt{14}$ \\
\bottomrule
\end{tabular}
\caption{For each point from the 5-point set $S_{+}$ in Fig.~\ref{fig:5-point_sets}, the distances to neighbors from Table~\ref{tab:distances_S-+} are ordered in each row.}
\label{tab:ordered_distances_S-+}
\end{table}

If we ignore labels of points, $S_{\pm}$ have identical Pointwise Distance Distribution ($\PDD$), which is the Simplexwise Distance Distributions ($\SDD$) in Definition~\ref{dfn:SDD}.
\smallskip

For easier visualization, the matrix below is obtained by lexicographically ordering the rows in Table~\ref{tab:ordered_distances_S-+}:
$$\PDD(S_{\pm})=\SDD(S_{\pm};1)=\left( \begin{array}{ccccc}
\sqrt{2} & \sqrt{6} & \sqrt{6} & \sqrt{14} \\
\sqrt{2} & \sqrt{6} & \sqrt{8} & \sqrt{14} \\
\sqrt{6} & \sqrt{6} & \sqrt{8} & \sqrt{14} \\
\sqrt{6} & \sqrt{6}  & \sqrt{14} & \sqrt{32}\\
\sqrt{6} & \sqrt{14}  & \sqrt{14} & \sqrt{32} 
\end{array} \right)$$
Now we show that $\SDD(S_-;2)\neq \SDD(S_+;2)$.
For $h=2$, the Simplexwise Distance Distribution $\SDD(C;h)$ consists of $\RDD(C;A)$ for 2-point subsets $A\subset C$.
Both sets $S_{\pm}$ have a single pair of points $(G_+,B-)$ and $(G+,B_+)$ at distance $\sqrt{2}$. 
Hence it suffices to show that the Relative Distance Distributions differ for this pair:
$$\RDD(S_-,\vect{G_+}{B_-})=\left[\sqrt{2},
\left( \begin{array}{ccc}
\sqrt{8} & \sqrt{14} & \sqrt{6} \\
\sqrt{6} & \cbox{yellow}{\sqrt{6}} & \cbox{yellow}{\sqrt{14}} \\
G_- & R_- & R_+ \\
\end{array} \right) \right]$$
$$\RDD(S_+,\vect{G_+}{B_+})=\left[\sqrt{2},
\left( \begin{array}{ccc}
\sqrt{8} & \sqrt{14} & \sqrt{6} \\
\sqrt{6} & \cbox{yellow}{\sqrt{14}} & \cbox{yellow}{\sqrt{6}} \\
G_- & R_- & R_+ \\
\end{array} \right) \right]$$

The last rows in the above $3\times 3$ matrices indicate a complementary point $q\in C-A$ for indexing columns of the $2\times 3$ matrices $R(C;A)$ in Definition~\ref{dfn:RDD}. 
The resulting $\RDD$s differ because any permutation of rows or columns of $R(S_+;\{G_+,B_+\})$ keeps the pair $\sqrt{6},\sqrt{6}$ in the same column but $R(S_-;\{G_+,B_+\})$ has no pair $\sqrt{6},\sqrt{6}$ in one column.
Hence $\SDD(S_-;2)\neq\SDD(S_+;2)$. 
%\bs
\end{exa}

%B2-----------------------------------------

\begin{figure*}[h!]
\centering
\includegraphics[width=\linewidth]{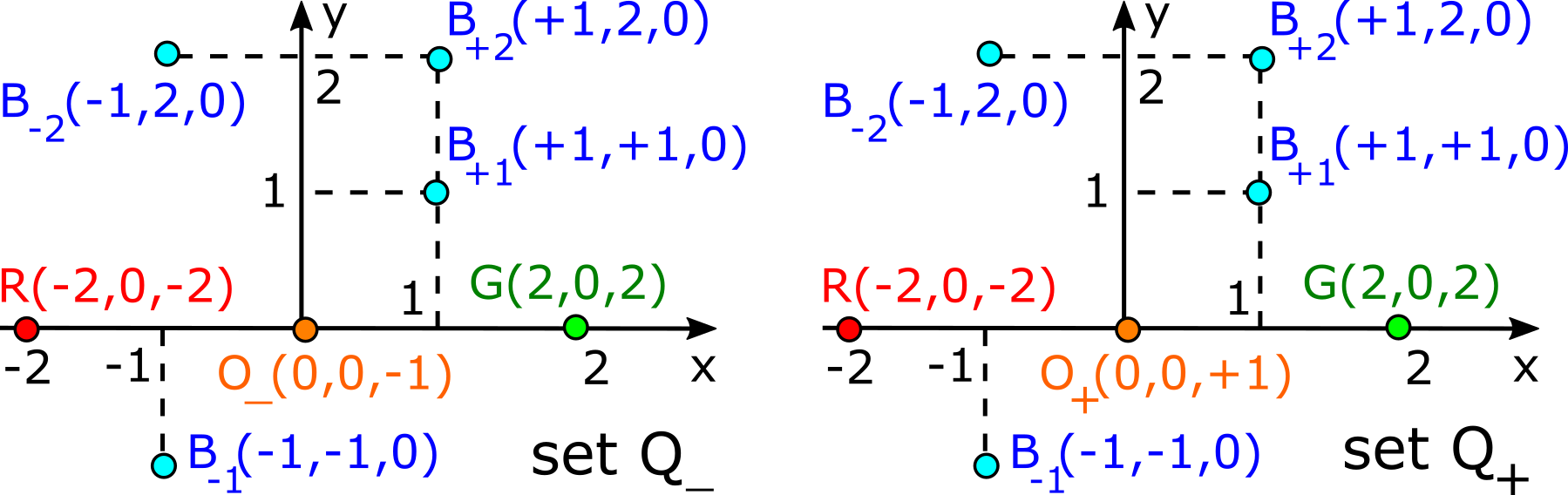}
\caption{\textbf{Left}: $(x,y)$-projection of the 7-point set $Q_-\subset\R^3$ of the red point $R=(-2,0,-2)$, green point $G=(2,0,2)$, four blue points $B_{\pm 1}=(\pm 1,\pm 1,0)$, $B_{\pm 2}=(\pm 1,2,0)$, orange point $O_-=(0,0,-1)$.
\textbf{Right}: to get $Q_+$ from $Q_-$, replace $O_-$ with $O_+=(0,0,+1)$.}
\label{fig:7-point_sets}
\end{figure*}

\begin{exa}[7-point sets]
\label{exa:7-point_sets}
The sets $Q_{\pm}$ in Fig.~\ref{fig:7-point_sets} taken from \cite[Figure S4(B)]{pozdnyakov2020incompleteness} have distances in Table~\ref{tab:distances_Q-+}.
Both sets have only two pairs of points at distance $\sqrt{6}$.
Hence it suffices to compare $\RDD$s for these pairs below.

$$R(Q_-;\vect{G}{B_{+1}})=
\left( \begin{array}{ccccc}
\sqrt{32} & \sqrt{14} & \sqrt{17} & 3 & \cbox{yellow}{\sqrt{13}}  \\
\sqrt{6} & \sqrt{8} & \sqrt{5} & 1 & \sqrt{3} \\
R & B_{-1} & B_{-2} & B_{+2} & O_-
\end{array} \right),$$
$$R(Q_-;\vect{R}{B_{-1}})=
\left( \begin{array}{ccccc}
\sqrt{32} & \sqrt{14} & 3 & \sqrt{17} & \cbox{yellow}{\sqrt{5}}  \\
\sqrt{14} & \sqrt{8} & 3 & \sqrt{13} & \sqrt{3} \\
G & B_{+1} & B_{-2} & B_{+2} & O_-
\end{array} \right)$$

The pair above has submatrices $\mat{3}{\sqrt{13}}{1}{\sqrt{3}}$ and $\mat{3}{\sqrt{5}}{3}{\sqrt{3}}$ but the pair below has no such submatrices.

$$R(Q_+;\vect{G}{B_{+1}})=
\left( \begin{array}{ccccc}
\sqrt{32} & \sqrt{14} & \sqrt{17} & 3 & \cbox{yellow}{\sqrt{5}}  \\
\sqrt{6} & \sqrt{8} & \sqrt{5} & 1 & \sqrt{3} \\
R & B_{-1} & B_{-2} & B_{+2} & O_+
\end{array} \right),$$
$$R(Q_+;\vect{R}{B_{-1}})=
\left( \begin{array}{ccccc}
\sqrt{32} & \sqrt{14} & 3 & \sqrt{17} & \cbox{yellow}{\sqrt{13}}  \\
\sqrt{14} & \sqrt{8} & 3 & \sqrt{13} & \sqrt{3} \\
G & B_{+1} & B_{-2} & B_{+2} & O_+
\end{array} \right)$$

The pair of $\RDD(Q_-;\vect{G}{B_{+1}})$ and $\RDD(Q_-;\vect{R}{B_{-1}})$ differs from the pair $\RDD(Q_+;\vect{G}{B_{+1}})$ and $\RDD(Q_+;\vect{R}{B_{-1}})$.
Hence $\SDD(Q_-;2)\neq\SDD(Q_+;2)$. 
%\bs
\end{exa}

\begin{table*}[h!]
\centering
\begin{tabular}{@{}l|ccccccc@{}}
\toprule
distances of \hl{ $Q_-$ } & $R$ & $G$ & $B_{-1}$ & $B_{+1}$ & $B_{-2}$ & $B_{+2}$ & \hl{ $O_{-}$ } \\
\midrule
$R=(-2,0,-2)$ & $0$ & $\sqrt{32}$ & $\sqrt{6}$ & $\sqrt{14}$ & $3$ & $\sqrt{17}$ & \hl{ $\sqrt{5}$ } \smallskip \\
$G=(+2,0,+2)$ & $\sqrt{32}$ & 0 & $\sqrt{14}$ & $\sqrt{6}$ & $\sqrt{17}$ & $3$ & \hl{ $\sqrt{13}$ } \smallskip \\
$B_{-1}=(-1,-1,0)$ & $\sqrt{6}$ & $\sqrt{14}$ & 0 & $\sqrt{8}$ & $3$ & $\sqrt{13}$ & $\sqrt{3}$ \smallskip \\
$B_{+1}=(+1,+1,0)$ & $\sqrt{14}$ & $\sqrt{6}$ & $\sqrt{8}$ & 0 & $\sqrt{5}$ & $1$ & $\sqrt{3}$ \smallskip \\
$B_{-2}=(-1,2,0)$ & $3$ & $\sqrt{17}$ & $3$ & $\sqrt{5}$ & 0 & 2 & $\sqrt{6}$ \smallskip \\
$B_{+2}=(+1,2,0)$ & $\sqrt{17}$ & $3$ & $\sqrt{13}$ & $1$ & 2 & 0 & $\sqrt{6}$ \smallskip \\
\hl{ $O_{-}=(0,0,-1)$ } & \hl{ $\sqrt{5}$ } & \hl{ $\sqrt{13}$ } & $\sqrt{3}$ & $\sqrt{3}$ & $\sqrt{6}$ & $\sqrt{6}$ & 0 \smallskip \\
\bottomrule
\end{tabular}
\smallskip

\begin{tabular}{@{}l|ccccccc@{}}
\toprule
distances of \hl{ $Q_+$ } & $R$ & $G$ & $B_{-1}$ & $B_{+1}$ & $B_{-2}$ & $B_{+2}$ & \hl{ $O_{+}$ } \\
\midrule
$R=(-2,0,-2)$ & $0$ & $\sqrt{32}$ & $\sqrt{6}$ & $\sqrt{14}$ & $3$ & $\sqrt{17}$ & \hl{ $\sqrt{13}$ } \smallskip \\
$G=(+2,0,+2)$ & $\sqrt{32}$ & 0 & $\sqrt{14}$ & $\sqrt{6}$ & $\sqrt{17}$ & $3$ & \hl{ $\sqrt{5}$ } \smallskip \\
$B_{-1}=(-1,-1,0)$ & $\sqrt{6}$ & $\sqrt{14}$ & 0 & $\sqrt{8}$ & $3$ & $\sqrt{13}$ & $\sqrt{3}$ \smallskip \\
$B_{+1}=(+1,+1,0)$ & $\sqrt{14}$ & $\sqrt{6}$ & $\sqrt{8}$ & 0 & $\sqrt{5}$ & $1$ & $\sqrt{3}$ \smallskip \\
$B_{-2}=(-1,2,0)$ & $3$ & $\sqrt{17}$ & $3$ & $\sqrt{5}$ & 0 & 2 & $\sqrt{6}$ \smallskip \\
$B_{+2}=(+1,2,0)$ & $\sqrt{17}$ & $3$ & $\sqrt{13}$ & $1$ & 2 & 0 & $\sqrt{6}$ \smallskip \\
\hl{ $O_{+}=(0,0,+1)$ } & \hl{ $\sqrt{13}$ } & \hl{ $\sqrt{5}$ } & $\sqrt{3}$ & $\sqrt{3}$ & $\sqrt{6}$ & $\sqrt{6}$ & 0 \\
\bottomrule
\end{tabular}
\caption{The matrices of distances between all points of the 7-point set $Q_{\mp}$ in Fig.~\ref{fig:7-point_sets} taken from \cite[Figure S4(B)]{pozdnyakov2020incompleteness}.}
\label{tab:distances_Q-+}
\end{table*}

%B3-----------------------------------------

\begin{exa}[6-point sets]
\label{exa:6-point_sets}
The sets $T_{\pm}$ in Fig.~\ref{fig:6-point_sets}, which was motivated by \cite[Figure S4(C)]{pozdnyakov2020incompleteness}, have the points $R,G,O_{\pm}$ from the sets $Q_{\pm}$ in Example~\ref{exa:7-point_sets} and three new points $C_1(x_1,y_1,0)$, $C_2(x_2,y_2,0)$, $C_3(x_3,y_3,0)$ such that $|RC_1|=|GC_2|$, $|RC_2|=|GC_3|$, $|RC_3|=|GC_1|$.

\begin{figure*}[h!]
\centering
\includegraphics[width=\linewidth]{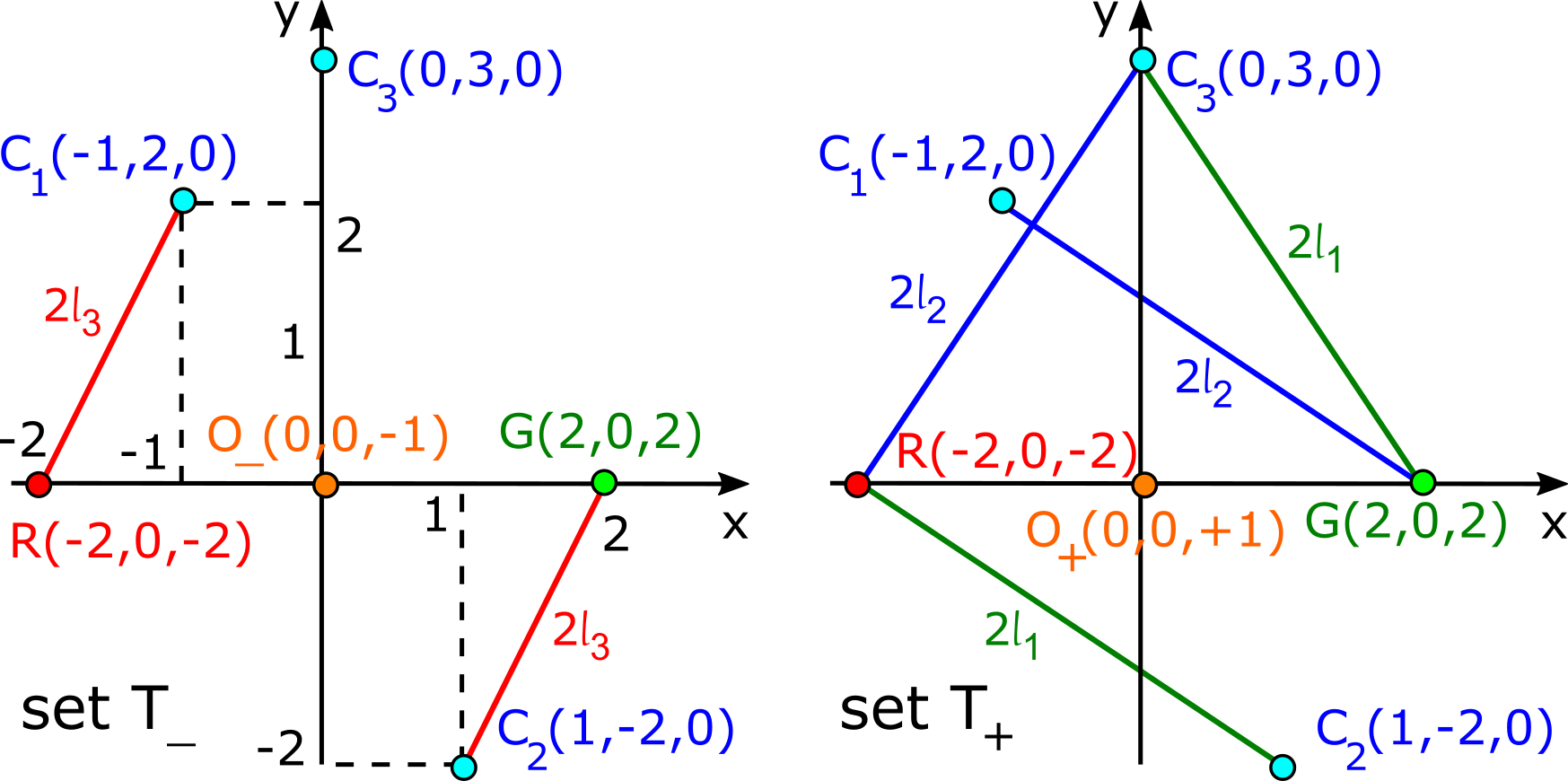}
\caption{\textbf{Left}: $(x,y)$-projection of the 6-point set $T_-\subset\R^3$ consisting of the red point $R=(-2,0,-2)$, green point $G=(2,0,2)$, three blue points $C_1=(x_1,y_1,0)$, $C_2=(x_2,y_2,0)$,  $C_3=(x_3,y_3,0)$ and orange point $O_-=(0,0,-1)$ so that $|RC_1|=2l_3=|GC_2|$, $|RC_2|=2l_1=|GC_3|$, $|RC_3|=2l_2=|GC_1|$.
\textbf{Right}: to get the set $T_+\subset\R^3$ from $T_-$, replace $O_-$ with $O_+=(0,0,+1)$.}
\label{fig:6-point_sets}
\end{figure*}

\begin{table*}[h!]
\centering
\begin{tabular}{@{}l|cccccc@{}}
\toprule
distances of \hl{ $T_-$ } & $R$ & $G$ & $C_{1}$ & $C_{2}$ & $C_{3}$ & \hl{ $O_{-}$ } \\
\midrule
$R=(-2,0,-2)$ & $0$ & $\sqrt{32}$ & $2\sqrt{l_3^2+1}$ & $2\sqrt{l_1^2+1}$ & $2\sqrt{l_2^2+1}$ & \hl{ $\sqrt{5}$ } \smallskip \\
$G=(+2,0,+2)$ & $\sqrt{32}$ & 0 & $2\sqrt{l_2^2+1}$ & $2\sqrt{l_3^2+1}$ & $2\sqrt{l_1^2+1}$  & \hl{ $\sqrt{13}$ } \smallskip \\
$C_{1}=(x_1,y_1,0)$ & $2\sqrt{l_3^2+1}$ & $2\sqrt{l_2^2+1}$ & 0 & $|C_1C_2|$ & $|C_3C_1|$ & $\sqrt{2l_2^2+2l_3^2-3}$ \smallskip \\
$C_{2}=(x_2,y_2,0)$ & $2\sqrt{l_1^2+1}$ & $2\sqrt{l_3^2+1}$ & $|C_1C_2|$ & 0 & $|C_2C_3|$ & $\sqrt{2l_3^2+2l_1^2-3}$ \smallskip \\
$C_{3}=(x_3,y_3,0)$ & $2\sqrt{l_2^2+1}$ & $2\sqrt{l_1^2+1}$ & $|C_3C_1|$ & $|C_2C_3|$ & 0 & $\sqrt{2l_1^2+2l_2^2-3}$ \smallskip  \\
\hl{ $O_{-}=(0,0,-1)$ } & \hl{ $\sqrt{5}$ } & \hl{ $\sqrt{13}$ } & $\sqrt{2l_2^2+2l_3^2-3}$ & $\sqrt{2l_3^2+2l_1^2-3}$ & $\sqrt{2l_1^2+2l_2^2-3}$ & 0 \smallskip \\
\bottomrule
\end{tabular}
\smallskip

\begin{tabular}{@{}l|cccccc@{}}
\toprule
distances of \hl{ $T_+$ } & $R$ & $G$ & $C_{1}$ & $C_{2}$ & $C_{3}$ & \hl{ $O_{+}$ } \\
\midrule
$R=(-2,0,-2)$ & $0$ & $\sqrt{32}$ & $2\sqrt{l_3^2+1}$ & $2\sqrt{l_1^2+1}$ & $2\sqrt{l_2^2+1}$ & \hl{ $\sqrt{13}$ } \smallskip \\
$G=(+2,0,+2)$ & $\sqrt{32}$ & 0 & $2\sqrt{l_2^2+1}$ & $2\sqrt{l_3^2+1}$ & $2\sqrt{l_1^2+1}$  & \hl{ $\sqrt{5}$ } \smallskip \\
$C_{1}=(x_1,y_1,0)$ & $2\sqrt{l_3^2+1}$ & $2\sqrt{l_2^2+1}$ & 0 & $|C_1C_2|$ & $|C_3C_1|$ & $\sqrt{2l_2^2+2l_3^2-3}$ \smallskip \\
$C_{2}=(x_2,y_2,0)$ & $2\sqrt{l_1^2+1}$ & $2\sqrt{l_3^2+1}$ & $|C_1C_2|$ & 0 & $|C_2C_3|$ & $\sqrt{2l_3^2+2l_1^2-3}$ \smallskip \\
$C_{3}=(x_3,y_3,0)$ & $2\sqrt{l_2^2+1}$ & $2\sqrt{l_1^2+1}$ & $|C_3C_1|$ & $|C_2C_3|$ & 0 & $\sqrt{2l_1^2+2l_2^2-3}$  \smallskip \\
\hl{ $O_{+}=(0,0,+1)$ } & \hl{ $\sqrt{13}$ } & \hl{ $\sqrt{5}$ } & $\sqrt{2l_2^2+2l_3^2-3}$ & $\sqrt{2l_3^2+2l_1^2-3}$ & $\sqrt{2l_1^2+2l_2^2-3}$ & 0 \smallskip \\
\bottomrule
\end{tabular}
\caption{The matrices of distances between all points of the 6-point set $T_{\mp}$ in Fig.~\ref{fig:6-point_sets} motivated by \cite[Figure S4(C)]{pozdnyakov2020incompleteness}.}
\label{tab:distances_T-+}
\end{table*}

Denote by $2l_1,2l_2,2l_3$ the lengths of these three pairs of line segments after their projection to the $xy$-plane so that
$$\leqno{(\ref{exa:6-point_sets}.1)}
\left\{\begin{array}{l}
(x_2+2)^2+y_2^2=|RC_2|^2-4=(2l_1)^2, \\
(x_3-2)^2+y_3^2=|GC_3|^2-4=(2l_1)^2;
\end{array} \right.$$ 
$$\leqno{(\ref{exa:6-point_sets}.2)}
\left\{\begin{array}{l}
(x_3+2)^2+y_3^2=|RC_3|^2-4=(2l_2)^2, \\
(x_1-2)^2+y_1^2=|GC_1|^2-4=(2l_2)^2;
\end{array} \right.$$ 
$$\leqno{(\ref{exa:6-point_sets}.3)}
\left\{\begin{array}{l}
(x_1+2)^2+y_1^2=|RC_1|^2-4=(2l_3)^2, \\
(x_2-2)^2+y_2^2=|GC_2|^2-4=(2l_3)^2.
\end{array} \right.$$ 

Comparing the first part of $(\ref{exa:6-point_sets}.1)$ with the second part side of $(\ref{exa:6-point_sets}.3)$, we get $(2l_1)^2-4x_2=(2l_3)^2+4x_2$, so $x_2=\dfrac{l_1^2-l_3^2}{2}$.
Similarly, 
$x_3=\dfrac{l_2^2-l_1^2}{2}$, $x_1=\dfrac{l_3^2-l_2^2}{2}$ so that $x_1+x_2+x_3=0$.
From the second part of $(\ref{exa:6-point_sets}.2)$, we get
$x_1^2-4x_1+4+y_1^2=4l_2^2$, so
$$|O_{\pm} C_1|^2=x_1^2+y_1^2+1=4l_2^2+4x_1-3=2l_2^2+2l_3^2-3,$$
$$\text{similarly }|O_{\pm} C_2|^2=2l_3^2+2l_1^2-3,\;
|O_{\pm} C_3|^2=2l_1^2+2l_2^2-3.$$
Then $|C_1C_2|^2=(x_1-x_2)^2+(y_1-y_2)^2=x_1^2+y_1^2$.
\smallskip

\begin{table*}[h!]
\centering
\begin{tabular}{@{}l|c|c|c@{}}
\toprule
\hl{ $T_-$ } pair & distance & common pairs in $\wSDD(T_{\pm};2)$
& pairs that differ in $\wSDD(T_{+};2)$ \smallskip \\
\midrule
$\{R,O_{-}\}$ & $\sqrt{5}$ & $(\sqrt{13},\sqrt{32})$ to $G$ &
$\begin{array}{l}
(2\sqrt{l_3^2+1},\sqrt{2l_2^2+2l_3^2-3}) \text{ to } C_1, \smallskip \\
(2\sqrt{l_1^2+1},\sqrt{2l_3^2+2l_1^2-3}) \text{ to } C_2, \smallskip \\
(2\sqrt{l_2^2+1},\sqrt{2l_1^2+2l_2^2-3}) \text{ to } C_3  \smallskip
\end{array}$ \\
\midrule
$\{G,O_{-}\}$ & $\sqrt{13}$ & $(\sqrt{5},\sqrt{32})$ to $R$ &
$\begin{array}{c}
(2\sqrt{l_2^2+1},\sqrt{2l_2^2+2l_3^2-3}) \text{ to } C_1, \smallskip \\
(2\sqrt{l_3^2+1},\sqrt{2l_3^2+2l_1^2-3}) \text{ to } C_2, \smallskip\\
(2\sqrt{l_1^2+1},\sqrt{2l_1^2+2l_2^2-3}) \text{ to } C_3 \smallskip 
\end{array}$ \\
\midrule
$\{R,C_{i+1}\}$ & $2\sqrt{l_{i}^2+1}$ &
$\begin{array}{l}
(2\sqrt{l_{i-1}^2+1},\sqrt{32}) \text{ to } G,\smallskip \\ 
(2\sqrt{l_{i+1}^2+1},|C_{i+1} C_{i-1}|) \text{ to } C_{i-1},\smallskip \\
(2\sqrt{l_{i-1}^2+1},|C_i C_{i+1}|) \text{ to } C_i \smallskip 
\end{array}$ &
$(\sqrt{5},\sqrt{2l_{i-1}^2+2l_i^2-3})$ to $O_-$ \\
\midrule
$\{G,C_{i-1}\}$ & $2\sqrt{l_i^2+1}$ &
$\begin{array}{l}
(2\sqrt{l_{i+1}^2+1},\sqrt{32}) \text{ to } R, \smallskip \\ 
(2\sqrt{l_{i-1}^2+1},|C_{i+1} C_{i-1}|) \text{ to } C_{i+1}, \smallskip \\
(2\sqrt{l_{i+1}^2+1},|C_{i-1} C_i|) \text{ to } C_i \smallskip 
\end{array}$ &
$(\sqrt{13},\sqrt{2l_i^2+2l_{i+1}^2-3})$ to $O_-$ \\
\bottomrule
\end{tabular}
\bigskip

\begin{tabular}{@{}l|c|c|c@{}}
\toprule
\hl{ $T_+$ } pair & distance & common pairs in $\wSDD(T_{\pm};2)$
& pairs that differ in $\wSDD(T_{-};2)$  \\
\midrule
$\{G,O_{+}\}$ & $\sqrt{5}$ & $(\sqrt{13},\sqrt{32})$ to $R$ &
$\begin{array}{c}
(2\sqrt{l_2^2+1},\sqrt{2l_2^2+2l_3^2-3}) \text{ to } C_1, \smallskip \\ 
(2\sqrt{l_3^2+1},\sqrt{2l_3^2+2l_1^2-3}) \text{ to } C_2, \smallskip \\ 
(2\sqrt{l_1^2+1},\sqrt{2l_1^2+2l_2^2-3}) \text{ to } C_3 \smallskip 
\end{array}$ \\
\midrule
$\{R,O_{+}\}$ & $\sqrt{13}$ & $(\sqrt{5},\sqrt{32})$ to $G$ &
$\begin{array}{c}
(2\sqrt{l_3^2+1},\sqrt{2l_2^2+2l_3^2-3}) \text{ to } C_1, \smallskip \\ 
(2\sqrt{l_1^2+1},\sqrt{2l_3^2+2l_1^2-3}) \text{ to } C_2, \smallskip \\ 
(2\sqrt{l_2^2+1},\sqrt{2l_1^2+2l_2^2-3}) \text{ to } C_3 \smallskip
\end{array}$ \\
\midrule
$\{R,C_{i+1}\}$ & $2\sqrt{l_i^2+1}$ &
$\begin{array}{l}
(2\sqrt{l_{i-1}^2+1},\sqrt{32}) \text{ to } G, \smallskip  \\ 
(2\sqrt{l_{i+1}^2+1},|C_{i+1} C_{i-1}|) \text{ to } C_{i-1}, \smallskip \\
(2\sqrt{l_{i-1}^2+1},|C_i C_{i+1}|) \text{ to } C_i \smallskip 
\end{array}$ &
$(\sqrt{13},\sqrt{2l_{i-1}^2+2l_i^2-3})$ to $O_+$ \\
\midrule
$\{G,C_{i-1}\}$ & $2\sqrt{l_i^2+1}$ &
$\begin{array}{l}
(2\sqrt{l_{i+1}^2+1},\sqrt{32}) \text{ to } R, \smallskip \\ 
(2\sqrt{l_{i-1}^2+1},|C_{i+1} C_{i-1}|) \text{ to } C_{i+1}, \smallskip \\
(2\sqrt{l_{i+1}^2+1},|C_{i-1} C_i|) \text{ to } C_i \smallskip 
\end{array}$ &
$(\sqrt{5},\sqrt{2l_i^2+2l_{i+1}^2-3})$ to $O_+$ \\
\midrule
\end{tabular}
\caption{For the sets $T_{\pm}$, the distributions $\wSDD(T_{\pm};2)$ can differ only by $\wRDD$s of the pairs $\{R,O_{\pm}\}, \{G,O_{\pm}\}, \{R,C_i\}, \{G,C_i\}$ shown above, where $i\in\{1,2,3\}$ is considered modulo 3 so that $1-1\equiv 3\pmod{3}$.
In rows of corresponding pairs of points, some pairs of distances are the same in both $\wSDD(T_{\pm};2)$, but other pairs can differ.
If parameters $l_1,l_2,l_3$ are pairwise distinct, the rows $\{R,O_-\},\{G,O_+\}$ include three different pairs of distances, so $\wSDD(T_{-};2)\neq \wSDD(T_{+};2)$, see Example~\ref{exa:6-point_sets}.}
\label{tab:distance_pairs_T-+}
\end{table*}

\begin{table*}[h!]
\centering
\begin{tabular}{@{}l|c|cccc|@{}}
\toprule
\hl{ $T_-$ } pair & distance &  distance to neighbor1 & distance to neighbor 2 & distance to neighbor 3 & distance to neighbor 4 \\
\midrule
$\{R,C_{1}\}$ & $3$ 
& $(\sqrt{2},\sqrt{17})$ to $C_3$ 
& \hl{ $(\sqrt{5},\sqrt{6})$ } to $O_{-}$ 
& $(\sqrt{17},\sqrt{20})$ to $C_2$ 
& $(\sqrt{17},\sqrt{32})$ to $G$ \smallskip \\
\midrule
$\{G,C_{2}\}$ & $3$ 
& \hl{ $(\sqrt{6},\sqrt{13})$ } to $O_{-}$ 
& $(\sqrt{17},\sqrt{20})$ to $C_1$ 
& $(\sqrt{17},\sqrt{26})$ to $C_3$ 
& $(\sqrt{17},\sqrt{32})$ to $R$  \smallskip \\
\bottomrule
\end{tabular}
\smallskip

\begin{tabular}{@{}l|c|cccc|@{}}
\toprule
\hl{ $T_+$ } pair & distance &  distance to neighbor1 & distance to neighbor 2 & distance to neighbor 3 & distance to neighbor 4 \\
\midrule
$\{R,C_{1}\}$ & $3$ 
& $(\sqrt{2},\sqrt{17})$ to $C_3$ 
& \hl{ $(\sqrt{6},\sqrt{13})$ } to $O_{+}$ 
& $(\sqrt{17},\sqrt{20})$ to $C_2$ 
& $(\sqrt{17},\sqrt{32})$ to $G$  \smallskip\\
\midrule
$\{G,C_{2}\}$ & $3$ 
& \hl{ $(\sqrt{5},\sqrt{6})$ } to $O_{+}$ 
& $(\sqrt{17},\sqrt{20})$ to $C_1$ 
& $(\sqrt{17},\sqrt{26})$ to $C_3$ 
& $(\sqrt{17},\sqrt{32})$ to $R$ \smallskip \\
\bottomrule
\end{tabular}
\caption{The above rows show that $\SDD(T_-;2)\neq\SDD(T_+;2)$ for the sets $T_{\pm}$ with $C_1=(-1,2,0)$, $C_2=(1,-2,0)$, $C_3=(0,3,0)$ so that $l_1=l_2=\frac{\sqrt{5}}{2}$, $l_3=\frac{\sqrt{13}}{2}$ in Table~\ref{tab:distance_pairs_T-+}.}
\label{tab:distance_pairs_T-+ex}
\end{table*}

Table~\ref{tab:distances_T-+} contains all pairwise distances between the points of $T_{\mp}$.
We show that $T_{\pm}$ differ by the simplified invariants $\wSDD(T_{\pm};2)$ below.
In each column of $R(C;A)$, we additionally allow any permutation of elements independent of other columns, so we could order each column (a pair of distances) lexicographically.
Denote the resulting simplification of $\RDD$ by $\wRDD$.  
Then $\wSDD(T_{\pm};2)$ have identical $\wRDD$s for the 2-point subsets $A$ from the list $\{R,G\}, \{O_{\pm},C_i\}, \{C_i,C_j\}$ for distinct $i,j=1,2,3$.
\smallskip

For example, both $\wRDD(T_{\pm};\{R,G\})$ start with the distance $|R-G|=\sqrt{32}$ and then include the same four pairs $(\sqrt{5},\sqrt{13})$, $(2\sqrt{l_i^2+1},2\sqrt{l_{i-1}^2+1})$ for $i\in\{1,2,3\}$ modulo $3$, which should be ordered and written lexicographically.
Hence it makes sense to compare $\wSDD(T_{\pm};2)$ only by the remaining $\wRDD(T_{\pm};A)$ for $A$ from the list $\{R,O_{\pm}\}, \{G,O_{\pm}\}, \{R,C_i\}, \{G,C_j\}$ in Table~\ref{tab:distance_pairs_T-+}.
\smallskip

Without loss of generality assume that $l_1\geq l_2\geq l_3$.
If all the lengths are distinct, then $l_1>l_2>l_3$.
Then the rows for $\{R,O_{-}\}$ and $\{G,O_{+}\}$ differ in Table~\ref{tab:distance_pairs_T-+} even after ordering each pair so that a smaller distance precedes a larger one, and after writing all pairs lexicographically. 
So $\wSDD(T_-;2)\neq \wSDD(T_+;2)$ unless two of $l_i$ are equal.
\smallskip

If (say) $l_1=l_2$, the lexicographically ordered rows of $\{R,O_{-}\}$ and $\{G,O_{+}\}$ coincide in $\wSDD(T_{\pm};2)$, similarly for the rows of $\{G,O_{-}\}$ and $\{R,O_{+}\}$.
Hence it suffices to compare only the six rows for the remaining pairs $\{R,C_i\},\{G,C_j\}$ in $\wSDD(T_{\pm};2)$.
\medskip

For $l_1=l_2$, we get $x_3=\dfrac{l_2^2-l_1^2}{2}=0$ and 
$x_1=-x_2=\dfrac{l_3^2-l_2^2}{2}$.
In equation $(\ref{exa:6-point_sets}.3)$ the equality 
$(x_1+2)^2+y_1^2=(x_2-2)^2+y_2^2$ with $x_1=-x_2$ implies that $y_1^2=y_2^2$.
The even more degenerate case $l_1=l_2=l_3$, means that $x_1=x_2=x_3=0$ and $y_1^2=y_2^2=y_3^2$, hence at least two of $C_1,C_2,C_3$ should coincide.
The above contradiction means that it remains to consider only the case $l_1=l_2>l_3$ when $x_1=-x_2\neq 0=x_3$ and $y_1=\pm y_2$, see Fig.~\ref{fig:6-point_sets}.
\medskip

If $y_1=y_2$, the sets $T_{\pm}$ are isometric by $(x,y,z)\mapsto(-x,y,-z)$. 
If $y_1=-y_2$ and $y_3=0$, the sets $T_{\pm}$ are isometric by $(x,y,z)\mapsto(-x,-y,-z)$.
If $y_1=-y_2$ and $y_3\neq 0$, then $C_1=(x_1,y_1,0)$, $C_2=(-x_1,-y_1,0)$, $C_3\neq(0,0,0)$.
Then among the six remaining rows, only the rows of $\{R,C_1\}$, $\{G,C_2\}$ have points at the distance $2\sqrt{l_3^2+1}$, see  Table~\ref{tab:distance_pairs_T-+} for $i=3$ considered modulo 3.
Then $i+1\equiv 1\pmod{3}$, $i-1\equiv 2\pmod{3}$, so $l_{i+1}=l_1=l_2=l_{i-1}$.
\medskip

Looking at the rows of $\{R,C_1\}$, $\{G,C_2\}$, the three common pairs in each of $\SDD(T_{\pm};2)$ include the same distance $2\sqrt{l_1^2+1}=2\sqrt{l_2^2+1}$ but differ by $|C_{i-1} C_i|=|C_2 C_3|\neq |C_3C_1|=|C_i C_{i+1}|$ as $C_1=\pm C_2$, $C_3\neq(0,0,0)$.
\medskip

This couple of different rows implies that $\SDD(T_{-};2)\neq \SDD(T_{+};2)$ due to the swapped distances $\sqrt{5},\sqrt{13}$ in their remaining pairs, see Table~\ref{tab:distance_pairs_T-+ex} for the sets $T_{\pm}$ in Fig.~\ref{fig:6-point_sets} with
%$C_1=(-1,2,0)$, $C_2=(1,-2,0)$, $C_3=(0,3,0)$ with
$l_1=l_2=\frac{\sqrt{13}}{2}$, $l_3=\frac{\sqrt{5}}{2}$. 
%\bs
\end{exa}

%5=====================
\section{Continuous and computable metrics on SDD}
\label{sec:metrics}

%The section the Lipschitz

The $m-h$ permutable columns of the matrix $R(C;A)$ in $\RDD$ from Definition~\ref{dfn:RDD} can be interpreted as $m-h$ unlabelled points in $\R^h$.
Since any isometry is bijective, the simplest metric respecting bijections is the bottleneck distance (also called the Wasserstein distance $W_{\infty}$). 

\begin{dfn}[bottleneck distance $W_{\infty}$]
\label{dfn:bottleneck}
For any vector $v=(v_1,\dots,v_n)\in \R^n$, the \emph{Minkowski} norm is $||v||_{\infty}=\max\limits_{i=1,\dots,n}|v_i|$.
For any vectors or matrices $N,N'$ of the same size, the \emph{Minkowski} distance is $L_{\infty}(N,N')=\max\limits_{i,j}|N_{ij}-N'_{ij}|$.
For clouds $C,C'\subset\R^n$ of $m$ unlabelled points, 
the \emph{bottleneck distance} $W_{\infty}(C,C')=\inf\limits_{g:C\to C'} \sup\limits_{p\in C}||p-g(p)||_{\infty}$ is minimized over all bijections $g:C\to C'$.
%\bs
\end{dfn}

\begin{lem}[the max metric $M_{\infty}$ on $\RDD$s]
\label{lem:RDD+metric}
For any $m$-point clouds and ordered $h$-point sequences $A\subset C$ and $A'\subset C'$, set $d(\xi)=\max\{L_{\infty}(\xi(D(A)),D(A')),W_{\infty}(\xi(R(C;A)),R(C';A'))\}$ for a permutation $\xi\in S_h$ on $h$ points. 
Then the max metric $M_{\infty}(\RDD(C;A),\RDD(C';A'))=\min\limits_{\xi\in S_h}d(\xi)$
satisfies all metric axioms on $\RDD$s from Definition~\ref{dfn:RDD} and can be computed in time $O(h!(h^2 +m^{1.5}\log^h m) )$. 
%\bs
\end{lem}
\begin{proof}[Proof of Lemma~\ref{lem:RDD+metric}]
The first metric axiom says that $\RDD(C;A),\RDD(C';A')$ are equivalent by Definition~\ref{dfn:RDD} if and only if $M_{\infty}(\RDD(C;A),\RDD(C';A'))=0$ or $d(\xi)=0$ for some permutation $\xi\in S_h$.
Then $d(\xi)=0$ is equivalent to $\xi(D(A))=D(A')$ and $\xi(R(C;A))=R(C';A')$ up to a permutation of columns due to the first axiom for $W_{\infty}$.
The last two conclusions mean that the Relative Distance Distributions $\RDD(C;A),\RDD(C';A')$ are equivalent by Definition~\ref{dfn:RDD}.
The symmetry axiom follows since any permutation $\xi$ is invertible.
To prove the triangle inequality 
$M_{\infty}(\RDD(C;A),\RDD(C';A'))+$ \\
$M_{\infty}(\RDD(C'';A''),\RDD(C';A'))\geq$ \\ 
$M_{\infty}(\RDD(C;A),\RDD(C'';A''))$, let $\xi,\xi'\in S_h$ be optimal permutations for the $M_\infty$ values in the left-hand side above. 
The triangle inequality for $L_\infty$ says that \\
$L_{\infty}(\xi(D(A)),D(A'))+$ \\
$L_{\infty}(\xi'(D(A'')),D(A'))\geq $ \\
$L_{\infty}(\xi(D(A)),\xi'(D(A'')))=$ \\
$L_{\infty}(\xi'^{-1}\xi(D(A)),D(A''))$, similarly for the bottleneck distance $W_\infty$ from Definition~\ref{dfn:bottleneck}.
Taking the maximum of $L_\infty,W_\infty$ preserves the triangle inequality. 
Then $M_{\infty}(\RDD(C;A),\RDD(C'' ;A''))=\min\limits_{\xi\in S_h}d(\xi)$ cannot be larger than $d(\xi'^{-1}\xi)$ for the composition of the permutations above, so the triangle inequality holds for $M_\infty$. 
\smallskip

For a fixed permutation $\xi\in S_h$, the distance $L_\infty(\xi(D(A)),D(A'))$ requires $O(h^2)$ time.
The bottleneck distance $W_{\infty}(\xi(R(C;A)),R(C';A'))$ on the $h\times(m-h)$ matrices $\xi(R(C;A))$ and $R(C';A')$ with permutable columns can be considered as the bottleneck distance on clouds of $(m-h)$ unlabelled points in $\R^h$, so $W_\infty(\xi(R(C;A)),R(C';A'))$ needs only $O(m^{1.5}\log^h m)$ time by
\cite[Theorem~6.5]{efrat2001geometry}.
The minimization over all permutations $\xi\in S_h$ gives the factor $h!$ in the final time.
\end{proof}

%We will use only $h=n$ for Euclidean space $\R^n$, so the factor $h!$ in Lemma~\ref{lem:RDD+metric} is practically small for $n=2,3$.
For $h=1$ and a 1-point subset $A\subset C$, the matrix $D(A)$ is empty, so $d(\xi)=W_{\infty}(\xi(R(C;A)),R(C';A'))$.
The metric $M_{\infty}$ on $\RDD$s will be used for intermediate costs to get metrics on unordered collections of $\RDD$s ($\SDD$s) by using the standard tools in Definitions~\ref{dfn:LAC} and~\ref{dfn:EMD} below.  

\begin{dfn}[Linear Assignment Cost LAC {\cite{fredman1987fibonacci}}]
\label{dfn:LAC}
For any $k\times k$ matrix of costs $c(i,j)\geq 0$, $i,j\in\{1,\dots,k\}$, the \emph{Linear Assignment Cost}   
$\LAC=\frac{1}{k}\min\limits_{g}\sum\limits_{i=1}^k c(i,g(i))$ is minimized for all bijections $g$ on the indices $1,\dots,k$.
%\bs
\end{dfn}

The normalization factor $\frac{1}{k}$ in $\LAC$ makes this metric better comparable with $\EMD$ whose weights sum up to 1.
  
\begin{dfn}[Earth Mover's Distance on distributions]
\label{dfn:EMD}
Let $B=\{B_1,\dots,B_k\}$ be a finite unordered set of objects with weights $w(B_i)$, $i=1,\dots,k$.
Consider another set $D=\{D_1,\dots,D_l\}$ with weights $w(D_j)$, $j=1,\dots,l$.
Assume that a distance between $B_i,D_j$ is measured by a metric $d(B_i,D_j)$.
A \emph{flow} from $B$ to $D$ is a $k\times l$ matrix  whose entry $f_{ij}\in[0,1]$ represents a partial \emph{flow} from an object $B_i$ to $D_j$.
The \emph{Earth Mover's Distance} \cite{rubner2000earth} is the minimum of
$\EMD(B,D)=\sum\limits_{i=1}^{k} \sum\limits_{j=1}^{l} f_{ij} d(B_i,D_j)$ over $f_{ij}\in[0,1]$ subject to 
%the conditions
$\sum\limits_{j=1}^{l} f_{ij}\leq w(B_i)$ for $i=1,\dots,k$, 
$\sum\limits_{i=1}^{k} f_{ij}\leq w(D_j)$ for $j=1,\dots,l$, and
$\sum\limits_{i=1}^{k}\sum\limits_{j=1}^{l} f_{ij}=1$.
%\bs
\end{dfn}

The first condition $\sum\limits_{j=1}^{l} f_{ij}\leq w(B_i)$ means that not more than the weight $w(B_i)$ of the object $B_i$ `flows' into all $D_j$ via the flows $f_{ij}$, $j=1,\dots,l$. 
The second condition $\sum\limits_{i=1}^{k} f_{ij}\leq w(D_j)$ means that all flows $f_{ij}$ from $B_i$ for $i=1,\dots,k$ `flow' to $D_j$ up to its weight $w(D_j)$.
The last condition
$\sum\limits_{i=1}^{k}\sum\limits_{j=1}^{l} f_{ij}=1$ forces all $B_i$ to collectively `flow' into all $D_j$.  
$\LAC$ \cite{fredman1987fibonacci} and $\EMD$ \cite{rubner2000earth} can be computed in a near cubic time in the sizes of given sets of objects. 
\smallskip

Theorem~\ref{thm:SDD_metrics}(b) extends the $O(m^{1.5}\log^n m)$ algorithm for fixed clouds of $m$ unlabelled points in \cite[Theorem~6.5]{efrat2001geometry} to the harder case of isometry classes but keeps the polynomial time in $m$ for a fixed dimension $n$.

\begin{thm}[time of metrics on $\SDD$s]
\label{thm:SDD_metrics}
For any $m$-point clouds $C,C'$ in their own metric spaces and $h\geq 1$, let the Simplexwise Distance Distributions $\SDD(C;h)$ and $\SDD(C';h)$ consist of $k=\binom{m}{h}$ $\RDD$s with equal weights $\frac{1}{k}$ without collapsing identical $\RDD$s.
\smallskip

\noindent
\textbf{(a)}
Using the $k\times k$ matrix of costs computed by the metric $M_{\infty}$ between $\RDD$s from $\SDD(C;h)$ and $\SDD(C';h)$, 
the Linear Assignment Cost $\LAC$ from Definition~\ref{dfn:LAC} satisfies all metric axioms on $\SDD$s and can be computed in time $O(h!(h^2 +m^{1.5}\log^h m)k^2 + k^3\log k)$.
\smallskip

\noindent
\textbf{(b)}
Let $\SDD(C;h)$ and $\SDD(C';h)$ have a maximum size $l\leq k$ after collapsing identical $\RDD$s. Then $\EMD$ from Definition~\ref{dfn:EMD} satisfies all metric axioms  on $\SDD$s and is computed in time $O(h!(h^2 +m^{1.5}\log^h m) l^2 +l^3\log l)$.
\end{thm}
\begin{proof}
The Linear Assignment Cost ($\LAC$) from Definition~\ref{dfn:LAC} is symmetric because any bijective matching can be reversed.
The triangle inequality for $\LAC$ follows from the triangle inequality for the metric $M_{\infty}$ in Lemma~\ref{lem:RDD+metric} by using a composition of bijections $\SDD(C;h)\to\SDD(C';h)\to\SDD(C'';h)$ matching all $\RDD$s similarly to the proof of Lemma~\ref{lem:RDD+metric}.
The first metric axiom for LAC means that $\LAC=0$ if and only if there is a bijection $g:\SDD(C;h)\to\SDD(C';h)$ so that all matched $\RDD$s are at distance $M_\infty=0$, so these $\RDD$s are equivalent (hence $\SDD$s are equal) due to the first axiom of $M_\infty=0$, which was proved in Lemma~\ref{lem:RDD+metric}.
\smallskip

The metric axioms for the Earth Mover's Distance ($\EMD$) are proved in the appendix of \cite{rubner2000earth} assuming the metric axioms for the underlying distance $d$, which is the metric $M_\infty$ from Lemma~\ref{lem:RDD+metric} in our case.
\smallskip

The time complexities for $\LAC$ and $\EMD$ follow from the time $O((h^2 +m^{1.5}\log^h m) h!)$ for $M_\infty$ in Lemma~\ref{lem:RDD+metric}, after multiplying by a quadratic factor for the size of cost matrices and adding a near cubic time 
\cite{fredman1987fibonacci,goldberg1987solving}.
\end{proof}

The Lipschitz continuity of $\SDD$ in Theorem~\ref{thm:SDD_continuity} needs
Lemma~\ref{lem:vector+inequality} follows from 
its partial case in Lemma~\ref{lem:pair+inequality}.

\begin{lem}
\label{lem:pair+inequality}
For any $a,b,c,d\in\R$, if $a\leq b$ and $c\leq d$ then 
$\max\{|a-c|,|b-d|\}\leq \max\{|a-d|,|b-c|\}$.
%\bs
\end{lem}
\begin{proof}
We consider several cases of the relative locations of the pairs $a\leq b$ and $c\leq d$ in the line $\R$.
\smallskip

\noindent  
Case $a\leq b\leq c\leq d$.
The required inequality follows from $\max\{c-a,d-b\}\leq d-a=\max\{d-a,c-b\}$.
\smallskip

\noindent  
Case $a\leq c\leq b\leq d$.
The required inequality follows from $\max\{c-a,d-b\}\leq d-a=\max\{d-a,b-c\}$.
\smallskip

\noindent  
Case $a\leq c\leq d\leq b$.
The inequality $\max\{c-a,b-d\}\leq\max\{d-a,b-c\}$ holds as $c-a\leq d-a$, $b-d\leq b-c$.
\smallskip

\noindent  
All other cases reduce to the cases above by the transformation $(a,b)\mapsto(-b,-a)$, $(c,d)\mapsto(-d,-c)$, which preserves the given condition and required conclusion.
\end{proof}

\begin{lem}[the metric $L_\infty$ respects ordering]
\label{lem:vector+inequality}
For any vector $v=(v_1,\dots,v_k)\in\R^k$, the vector $\vec v\in\R^k$ is obtained from $v$ by writing all coordinates in increasing order.
Then $|\vec u-\vec v|_\infty \leq |u-v|_\infty$ for any vectors $u,v\in\R^k$.
%\bs
\end{lem}
\begin{proof}
Since $|u-v|_\infty=\max\limits_{i=1,\dots,k}|u_i-v_i|$, the metric $L_\infty$ is preserved under any permutation $\xi\in S_k$ applied simultaneously to the coordinates of both $u,v$.
Hence, without loss of generality, we can assume that all coordinates of one vector $u$ are already in increasing order, so $u=\vec u$.
For any pair of successive coordinates $u_i\leq u_{i+1}$, let the corresponding pair in $v$ be in the opposite order $v_i>v_{i+1}$.
\smallskip

By Lemma~\ref{lem:pair+inequality} the swap $v_i\lra v_{i+1}$ does not increase the $L_{\infty}$ distance between $(u_i,u_{i+1})$ and $(v_i,v_{i+1})$, hence between $u,v\in\R^k$.
Applying such swaps puts all coordinates of $v$ in increasing order without increasing $L_{\infty}$.   
\end{proof}

Theorem~\ref{thm:SDD_continuity} substantially generalizes the fact that perturbing two points in their $\ep$-neighborhoods changes the Euclidean distance between these points by at most $2\ep$. 
%\smallskip

\begin{thm}[Lipschitz continuity of $\SDD$s]
\label{thm:SDD_continuity}
In any metric space, let $C'$ be obtained from a cloud $C$ by perturbing each point within its $\ep$-neighborhood.
For any $h\geq 1$, $\SDD(C;h)$ changes by at most $2\ep$ in the $\LAC$ and $\EMD$ metrics.
The lower bound holds: $\EMD\big(\SDD(C;h),\SDD(C';h)\big)\geq|\SDM(C;h,1)-\SDM(C';h,1)|_\infty$.
%\bs
\end{thm}
\begin{proof}
Order all points of the given clouds $C,C'$ so that every point $p_i\in C$ has the same index as its perturbation $p'_i\in C'$.
In the given metric space, the distance $d(p_i,p_j)$ between any points in $C$ changes under perturbation by at most $2\ep$ so that $|d(p_i,p_j)-d(p'_i,p'_j)|\leq 2\ep$.
This upper bound $2\ep$ remains for the max metric $M_\infty$ from Lemma~\ref{lem:RDD+metric}, also for the LAC and EMD metrics due to the total weight 1 of all costs in Definitions~\ref{dfn:LAC} and~\ref{dfn:EMD}, respectively. 
\smallskip

Lemma~\ref{lem:vector+inequality} implies that re-writing coordinates of a vector in increasing order cannot increase the metric $L_{\infty}$, hence
$L_{\infty}(\xi(D(A)),D(A'))\geq |\SDV(A)-\SDV(A')|_\infty$
for any permutation $\xi\in S_h$ of indices $1,\dots,h$.
\smallskip

The bottleneck distance $W_{\infty}(\xi(R(C;A)),R(C';A'))$ is the maximum of the metric $L_{\infty}$ computed between corresponding
column vectors of the $h\times(m-h)$ matrices $\xi(R(C;A))$ and $R(C';A')$.
Let $d=(d_1,\dots,d_h)$ and $d'=(d'_1,\dots,d'_h)$ be two such columns.
The triangle inequalities imply that $|d-d'|_\infty=\max\limits_{i=1,\dots,h}|d_i-d'_i|\geq \frac{1}{h}\sum\limits_{i=1}^h |d_i-d'_i|\geq |\frac{1}{h}\sum\limits_{i=1}^h d_i - \frac{1}{h}\sum\limits_{i=1}^h d'_i|$.
Hence taking averages of all vector coordinates cannot increase the metric $L_\infty$.
Then $W_{\infty}(\xi(R(C;A)),R(C';A'))$ has the lower bound equal to the metric $L_\infty$ between the vectors of $m-h$ column averages in the matrices $\xi(R(C;A))$ and $R(C';A')$.
Applying Lemma~\ref{lem:vector+inequality} to these vectors in $\R^{m-h}$ implies that
$W_{\infty}(\xi(R(C;A)),R(C';A'))\geq |\vec R(C;A)-\vec R(C';A')|_\infty$
\smallskip

Taking the maximum of the metrics $L_\infty$ and $W_\infty$ considered above, we get the lower bound in terms of the Average Distance Distribution from Definition~\ref{dfn:SDM}:
$d(\xi)=\max\{L_{\infty}(\xi(D(A)),D(A')),W_{\infty}(\xi(R(C;A)),R(C';A'))\}\geq |\ADD(C;A)-\ADD(C';A')|_\infty$
\smallskip

Since the above argument holds for any permutation $\xi\in S_h$, we get 
$M_{\infty}(\RDD(C;A),\RDD(C';A'))=\min d(\xi)\geq |\ADD(C;A)-\ADD(C';A')|_\infty$.
\smallskip

Both $\SDD(C;h)$ and $\ASD(C;h)$ are unordered collections of $\binom{m}{h}$ $\RDD$ and vectors, respectively.
If we use an optimal flow matrix $f_{ij}$ for $\EMD\big(\SDD(C;h),\SDD(C';\big)$ from Definition~\ref{dfn:EMD} to compute $\EMD$ on $\ASD$ vectors, we get an upper bound for $\EMD\big(\ASD(C;h),\ASD(C';h)\big)$, which can be potentially smaller (for anoth flow matrix) but not larger, so
$\EMD\big(\SDD(C;h),\SDD(C';h)\big)\geq \EMD\big(\ASD(C;h),\ASD(C';h)\big)$.
%\smallskip
Considering $\ASD(S;h)$ as a weighted distribution of vectors, $\SDM(C;h;1)$ is its centroid from section~3 in \cite[section~3]{cohen1997earth}.
%``The {E}arth {M}over's {D}istance: lower bounds and invariance under translation'' by S.~Cohen and L.~Guibas, tech. report at Stanford University (1997), https://apps.dtic.mil/sti/pdfs/ADA358270.pdf.
The lower bound $\EMD\big(\ASD(C;h),\ASD(C';h)\big)\geq|\SDM(C;h,1)-\SDM(C';h,1)|_\infty$ follows from \cite[Theorem~1]{cohen1997earth}.
\end{proof}

%The final lower bound for $T,K$ in Fig.~\ref{fig:4-point_clouds} is $1.55>0.26$.

%6=================
\section{Measured Simplexwise Distribution (MSD)}
\label{sec:MSD}

This section adapts $\SDD$ for metric-measure spaces.

\begin{dfn}[metric-measure space]
\label{dfn:mm-space}
According to Gromov \cite[section~$3\frac{1}{2}.1$]{gromov1999metric}, a \emph{metric-measure space} $(X,d_X,\mu_X)$ is a compact space $X$ with a metric $d_X$ and a Borel measure $\mu_X$ such that $\mu_X(X)<+\infty$.
An \emph{isomorphism} between metric-measure spaces is an isometry $f:X\to Y$ that respects the measures in the sense that $\mu_Y(U)=\mu_X(f^{-1}(X))$ for any subset $U\subset Y$.
\end{dfn}

Dividing $\mu_X(U)$ by the full measure $\mu_X(X)<+\infty$ for any $U\subset X$, we can assume that $\mu_X(X)=1$, so $\mu_X$ is a probability measure.
Any metric space $X$ of $m$ points can be considered a metric-measure space with the uniform measure $\mu_X(p)=\dfrac{1}{m}$ for all points $p\in X$.
\smallskip

On two points $0,1$ in Euclidean line $\R$, the mm-spaces $X=(\{0,1\},1,\{\frac{1}{2},\frac{1}{2}\})$ and $Y=(\{0,1\},1,\{\frac{1}{3},\frac{2}{3}\})$ are isometric but not isomorphic because of different weights.
\smallskip

%Problem~\ref{pro:isometry} becomes much harder if we replace isometries between metric spaces with isomorphisms between mm-spaces because all known isometry invariants should be further refined to distinguish up to isomorphism.

Definition~\ref{dfn:MSD} extends the local distribution of distances from \cite[Definition~5.5]{memoli2011gromov} to higher orders $h>1$.

\begin{dfn}[Measured Simplexwise Distribution $\MSD$]
\label{dfn:MSD}
Let $(X,d_X,\mu_X)$ be any \emph{metric-measure space}.
For any basis sequence $A=(p_1,\dots,p_h)\in X^h$ of $h\geq 1$ ordered points, write the triangular distance matrix $D(A)$ from Definition~\ref{dfn:RDD} row-by-row as the \emph{Vector of Interpoint Distances} $\VID(A)\in\R_+^{h(h-1)/2}$ so that $\VID_k=d_X(p_i,p_j)$ for $k=h(i-1)+j-1$, $1\leq i<j\leq h$.
For a vector $\vec d=(d_1,\dots,d_h)\in\R_+^h$ of distance thresholds, the \emph{Vector of Sequence-based Measures} $\VSM(A;\vec d)\in\R_+^h$ consists of $h$ values $\mu_X(\{q\in X\vl d_X(q,p_i)\leq d_i\})$ for $i=1,\dots,h$.  
The \emph{Measured Simplexwise Distribution} of order $h\geq 1$ is the function $\MSD[X;h]: X^h\times\R_+^h\to\R_+^{h(h+1)/2}$ mapping any $A\in X^h$ and $\vec d\in\R_+^h$ to the pair $[\VID(A),\VSM(A;\vec d)]$ considered as a concatenated vector in $\R_+^{h(h+1)/2}$.
\end{dfn}

For $h=1$, the vector $\VID(A)$ is empty and the Measured Simplexwise Distribution of order $h=1$ coincides with the local distribution of distances \cite[Definition~5.5]{memoli2011gromov} $\MSD[X;1]: X\times\R_+\to\R_+$ mapping any point $p\in X$ and a threshold $d\in\R_+$ to %the measure value 
$\mu_X(\{q\in X \vl d_X(q,p)\leq d\})$.
\smallskip

Any permutation $\xi$ on indices $1,\dots,h$ naturally permutes the components of $\MSD[X;h]$.
If $X$ consists of $m$ points, $\MSD[X;h]$ reduces to the finite collection
 of $\binom{m}{h}$ vectors $\VID(A)$ paired with fields $\VSM(A;\vec d):\R_+^h\to\R_+^{h}$ only for unordered $h$-point subsets $A\subset X$, which can be refined to a stronger invariant analog of $\SDD$ below.  
%\smallskip

\begin{figure*}[h!]
\centering
\includegraphics[width=\linewidth]{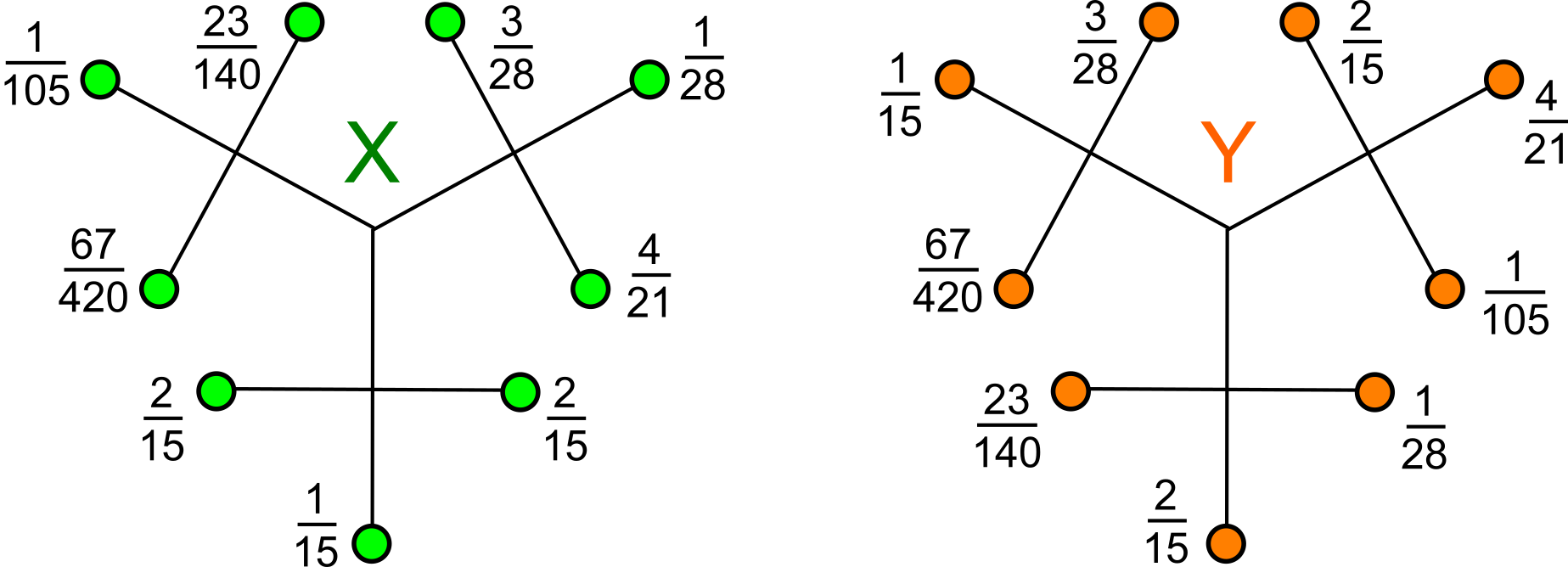}
\caption{Non-isomorphic mm-spaces $X,Y$ from \cite[Fig.~8]{memoli2011gromov} have equal local distributions of distances but are distinguished by the new Weighted Simplexwise Distribution of order 1 and the Measured Simplexwise Distributions of order 2, see details in Example~\ref{exa:9-point_trees}.
%All edges in trees have length $\frac{1}{2}$. Weights of points are indicated.
}
\label{fig:9-point_trees}
\end{figure*}

\begin{dfn}[Weighted Simplexwise Distribution $\WSD$]
\label{dfn:WSD}
Let $X$ be a finite mm-space whose any point $p$ has a weight $w_(p)$.
For $h\geq 1$ and a sequence $A=(p_1,\dots,p_h)\in X^h$ in Definition~\ref{dfn:RDD}, endow any distance $d(p,q)$ in $D(A)$ with the unordered pair $w(p),w(q)$ of weights.
For every point $q\in X-A$, put the weight $w(q)$ in the extra $(h+1)$-st row of the matrix $M(X;A)$ whose columns are indexed by unordered $q\in X-A$.
If $h=1$ and $A=p_1$, set $D(A)=w(p_1)$.
The \emph{Weighted Distance Distribution} $\WDD(X;A)$ is the equivalence class of the pair $[D(A);M(X;A)]$ up to permutations $\xi\in S_h$ acting on $A$. 
The \emph{Weighted Simplexwise Distribution} $\WSD(X)$ is the unordered collection of $\WDD(X;A)$ for all subsets $A\subset X$ of unordered $h$ points.
\end{dfn}

For finite mm-spaces, a metric on $\WDD$s can be defined similar to $M_\infty$ from Lemma~\ref{lem:RDD+metric} by combining the weights and distances.
Then $\LAC$ and $\EMD$ from Definitions~\ref{dfn:LAC} and~\ref{dfn:EMD} can be computed as in Theorem~\ref{thm:SDD_metrics}.

\begin{comment}
\begin{lem}[invariance of $\MSD$]
\label{lem:MSD_invariant}
The Measured Simplexwise Distribution $\MSD$ %[X;h]$ from Definition~\ref{dfn:MSD} 
is invariant under isomorphism.
% invariant of any mm-space $X$. 
\end{lem}
\begin{proof}
It easily follows from Definitions~\ref{dfn:mm-space} and~\ref{dfn:MSD}.
\end{proof}
\end{comment}

\begin{exa}[the strength of $\WSD$]
\label{exa:9-point_trees}
Fig.~\ref{fig:9-point_trees} shows mm-spaces $X,Y$ on 9 points visualised as trees \cite[Fig.~8]{memoli2011gromov}.
All edges have length $\frac{1}{2}$ and induce the shortest-path metrics $d_X,d_Y$.
The sum of weights in every small branch of 3 nodes is $\frac{1}{3}$.
These mm-spaces $X,Y$ have all inter-point distances only 1 and 2, and equal local distributions of distances $\MSD[X;1]=\MSD[Y;1]$ by \cite[Example~5.6]{memoli2011gromov}.
\smallskip

Indeed, both $\MSD$s can be considered the same set of 9 piecewise constant functions $\mu(p)$ taking values $w(p)$, $\frac{1}{3}$, and $1$ on the intervals $[0,1)$, $[1,2)$, $[2,+\infty)$, respectively.
\smallskip

However, $\WSD$s have more pointwise data: $\WSD[X;1]$ has $A(D)=w(p)=\frac{23}{140}$ and the following $2\times 8$ matrix 
$$M(X;p)=\left(\begin{array}{cccccccc}
1 & 1 & 2 & 2 & 2  & 2 & 2 & 2 \\
\frac{1}{105} & \frac{67}{240} & \frac{2}{15} & \frac{1}{15} & 
\frac{2}{15} & \frac{4}{21} & \frac{1}{28} &  \frac{3}{28}
\end{array}\right),$$ but 
$\WSD[Y;1]$ has another matrix for $w(p)=\frac{23}{140}$.
$$M(Y;p)=\left(\begin{array}{cccccccc}
1 & 1 & 2 & 2 & 2  & 2 & 2 & 2 \\
\frac{2}{15} & \frac{1}{28} & \frac{1}{105} & \frac{4}{21} & 
\frac{2}{15} & \frac{3}{28} & \frac{1}{15} & \frac{67}{420} 
\end{array}\right).$$
The matrices with freely permutable columns are different, so
 $X,Y$ are distinguished by $\WSD$ for $h=1$
\smallskip 

Also, $\MSD[X;2]\neq\MSD[Y;2]$ because, for any basis sequence $A=(p,q)\in X^2$, we have $\VSM[X;2](A;d_1,d_2)=(w(p),w(q))$ for $d_1,d_2<1$ since all other points have minimum distance $1$ from $p,q$, similarly for $Y$.
The unique points $p,q$ of weights $w(p)=\dfrac{23}{140}$ and $w(q)=\dfrac{67}{420}$ have different distances $d_X(p,q)=1$ and $d_Y(p,q)=2$.
Then $\MSD[X;2]\neq\MSD[Y;2]$ differ by the uniquely identifiable fields mapping $[0,1)^2$ to the constant vector $(w(p),w(q))$ with $\VID_X(A)=1\neq 2=\VID_Y(A)$.
\end{exa}

We conjecture that any mm-spaces $X,Y$ can be distinguished up to isomorphism by Measured Simplexwise Distributions for a high enough $h$ depending on $X,Y$.
\smallskip

Future updates of this paper will include continuous metrics between Measured Simplexwise Distributions on mm-spaces.
We are open to new ideas and collaboration.
\smallskip

New invariants in Definitions~\ref{dfn:SDD},~\ref{dfn:MSD} and main Theorems~\ref{thm:SDD_time},~\ref{thm:SDD_metrics},~\ref{thm:SDD_continuity} 
essentially contribute to the new area of \emph{Geometric Data Science} aiming to resolve all data challenges whose bottlenecks are analogues of Problem~\ref{pro:isometry}.
\smallskip

The earlier work has studied the following important cases of Problem~\ref{pro:isometry}: 1-periodic discrete series \cite{anosova2022density,anosova2023density,kurlin2022computable},
2D lattices \cite{kurlin2022mathematics,bright2023geographic}, 3D lattices \cite{mosca2020voronoi,bright2021welcome,kurlin2022exactly,kurlin2022complete}, periodic point sets in $\R^3$ \cite{smith2022practical,edelsbrunner2021density} and in higher dimensions \cite{anosova2021introduction,anosova2021isometry,anosova2022algorithms}.
\smallskip

The applications of Geometric Data Science to crystalline materials  \cite{ropers2022fast,balasingham2022compact,vriza2022molecular,zhu2022analogy}
led to the \emph{Crystal Isometry Principle} \cite{widdowson2022average,widdowson2021pointwise,widdowson2022resolving}
extending Mendeleev's table of chemical elements to the \emph{Crystal Isometry Space} of all periodic crystals continuously parametrised by complete invariants.
\smallskip

This work was supported by the Royal Academy of Engineering fellowship ``Data science for next generation engineering of solid crystalline materials'' (2021-2023, IF2122/186) and the EPSRC grants ``Application-driven Topological Data Analysis'' (2018-2023,  EP/R018472/1) and ``Inverse design of periodic crystals'' (2022-2024, EP/X018474/1).

The author thanks all members of the Data Science Theory and Applications group in the Materials Innovation Factory (Liverpool, UK), especially Daniel Widdowson, Matthew Bright, Yury Elkin, Olga Anosova, also Justin Solomon (MIT), Steven Gortler (Harvard),  Nadav Dym (Technion) for fruitful discussions, and any reviewers for their valuable time and helpful suggestions.

{\small
\bibliographystyle{ieee_fullname}
\bibliography{SDD}
}

\end{document}